\documentclass[10pt,a4wide]{amsart}
\usepackage[english]{babel}
\usepackage[T1]{fontenc}
\usepackage[latin1]{inputenc}
\usepackage{amssymb}
 \theoremstyle{plain}    
 \newtheorem{thm}{Theorem}[section]
 \theoremstyle{definition}
 \newtheorem{defn}[thm]{Definition}
 \theoremstyle{definition}
  \newtheorem{example}[thm]{Example}
 \theoremstyle{plain}    
 \newtheorem{prop}[thm]{Proposition} 
 \theoremstyle{plain}    
 \newtheorem{cor}[thm]{Corollary} 
 \theoremstyle{remark}    
 \newtheorem{notation}[thm]{Notation} 

\usepackage{latexsym,amsfonts,xypic}
\usepackage{pstricks,pst-node,pst-tree}
\xyoption{all}
\newcommand{\bc}{\mathbb{C}}

\newcommand{\ba}{\mathbb{A}}
\newcommand{\bp}{\mathbb{P}}

\theoremstyle{definition}
\newtheorem{enavant}[thm]{}

\begin{document}

\title{Danielewski-Fieseler surfaces }

\author{Adrien Dubouloz}
\address{Adrien Dubouloz \newline\indent Institut Fourier,    Laboratoire de Math\'ematiques
 \newline\indent  UMR5582 (UJF-CNRS) \newline\indent  BP 74 St Martin d'H\`eres
Cedex, FRANCE}
\email{ adrien.dubouloz@ujf-grenoble.fr}

\begin{abstract}
We study a class of normal affine surfaces with additive group actions
which contains in particular the Danielewski surfaces in $\ba^{3}$
given by the equations $x^{n}z=P\left(y\right)$, where $P$ is a
nonconstant polynomial with simple roots. We call them Danielewski-Fieseler
Surfaces. We reinterpret a construction of Fieseler \cite{Fie94}
to show that these surfaces appear as the total spaces of certain
torsors under a line bundle over a curve with an $r$-fold point.
We classify Danielewski-Fieseler surfaces through labelled rooted
trees attached to such a surface in a canonical way. Finally, we characterize
those surfaces which have a trivial Makar-Limanov invariant in terms
of the associated trees. 
\end{abstract}
\maketitle

\section*{Introduction}

The surfaces $S_{j}=\left\{ x^{j}z=y^{2}+1\right\} $, $j\geq1$,
in $\bc^{3}=\textrm{Spec}\left(\bc\left[x,y,z\right]\right)$ admit
free actions of the additive group $\mathbb{G}_{a,\bc}$ induced by
the locally nilpotent derivations $\partial_{j}=x^{j}\partial_{y}+2y\partial_{z}$
of $\bc\left[x,y,z\right]$ respectively. Danielewski \cite{Dan89}
observed that the associated quotient morphism $\pi:S_{j}\rightarrow S_{j}/\!/\mathbb{G}_{a,\bc}\simeq\textrm{Spec}\left(\bc\left[x\right]\right)$
is an $\ba^{1}$-fibration which factors throught an $\ba^{1}$-bundle
$\rho:S_{j}\rightarrow X$ over the affine line with a double origin
$\delta:X\rightarrow\textrm{Spec}\left(\bc\left[x\right]\right)$.
Indeed, $S_{j}$ is obtained as the gluing of two copies $\textrm{Spec}\left(\bc\left[x\right]\left[T_{i}\right]\right)$
of $\bc^{2}$, $i=1,2$, over $\bc^{*}\times\bc$ by means of the
$\bc\left[x\right]_{x}$-algebras isomorphisms $\bc\left[x\right]_{x}\left[T_{1}\right]\stackrel{\sim}{\rightarrow}\bc\left[x\right]_{x}\left[T_{2}\right]$,
$T_{1}\mapsto2x^{-j}+T_{2}$. This interpretation was further generalized
by Fieseler \cite{Fie94} to describe certain invariant neighbourhoods
of the fibers of a quotient $\ba^{1}$-fibration $\pi:S\rightarrow Z$
associated with a nontrivial action of $\mathbb{G}_{a,\bc}$ on a
normal affine surface $S$. More precisely, he established that if
a fiber $\pi^{-1}\left(z_{0}\right)$ is reduced then the induced
morphism $p_{2}:S\times_{Z}\textrm{Spec}\left(\mathcal{O}_{Z,z_{0}}\right)\rightarrow\textrm{Spec}\left(\mathcal{O}_{Z,z_{0}}\right)$,
where $\mathcal{O}_{Z,z_{0}}$ denotes the local ring of $z_{0}\in Z$,
factors through an $\ba^{1}$-bundle $\rho:S\rightarrow X$ over the
scheme $X$ obtained from $\textrm{Spec}\left(\mathcal{O}_{Z,z_{0}}\right)$
by replacing $z_{0}$ by as many points as there are connected components
in $\pi^{-1}\left(z_{0}\right)$. More generally, given a field $k$
of caracteristic zero and a pair $\left(X=\textrm{Spec}\left(A\right),x_{0}=\textrm{div}\left(x\right)\right)$,
where $A$ is either discrete valuation ring with uniformizing parameter
$x$ and residue field $k$ or of a polynomial ring $A=k\left[x\right]$,
the same description holds for the following class of surfaces. 

\begin{defn}
A \emph{Danielewski-Fieseler surface} \emph{with base} $\left(A,x\right)$
(A DFS, for short) is an integral affine $X$-scheme $\pi:S\rightarrow X$
of finite type, restricting to a trivial line bundle over $X_{*}=X\setminus\left\{ x_{0}\right\} $,
and such that $\pi^{-1}\left(x_{0}\right)$ is nonempty and reduced,
consisting of a disjoint union of curves isomorphic to affine line
$\ba_{k}^{1}$. 
\end{defn}
In this paper, we give a combinatorial description of DFS's in terms
of $\left(A,x\right)$\emph{-labelled rooted trees}, that is, pairs
$\gamma=\left(\Gamma,\sigma\right)$ consisting of a rooted tree $\Gamma$
and a cochain $\sigma\in A^{n}$, indexed by the terminal elements
of $\Gamma$, satisfying certain conditions with respect to the geometry
of $\Gamma$ (see (\ref{def:W_Rooted_Tree}) below). Then, as an application,
we characterize the DFS's with base $\left(k\left[x\right],x\right)$
which have a trivial Makar-Limanov invariant \cite{KML97}. 

The paper is divided as follows. In section $1$, we collect some
preliminary results on labelled rooted trees. In section $2$, we
reinterpret the 'cocycle construction' of Fieseler \cite{Fie94} to
describe DFS's as torsors under certain line bundles on a nonseparated
scheme. In section $3$, we associate to every labelled tree $\gamma$,
a DFS $\pi_{\gamma}:S\left(\gamma\right)\rightarrow X$ which comes
equipped with a canonical birational $X$-morphism $\psi_{\gamma}:S\left(\gamma\right)\rightarrow\ba_{X}^{1}$.
For instance, given a polynomial $P\in k\left[y\right]$ with simple
roots $y_{1},\ldots,y_{n}\in k$, the Danielewski-Fieseler surface
with base $\left(k\left[x\right],x\right)$ \[
\pi:S_{P,m}=\textrm{Spec}\left(k\left[x,y,z\right]/\left(x^{m}z-P\left(y\right)\right)\right)\longrightarrow X=\textrm{Spec}\left(k\left[x\right]\right),\]
 equipped with the morphism $pr_{y}:S_{P,m}\rightarrow\ba_{X}^{1}=\textrm{Spec}\left(k\left[x\right]\left[y\right]\right)$
corresponds to the following labelled tree $\gamma_{P,m}=\left(\Gamma_{m,n},\sigma\right)$.

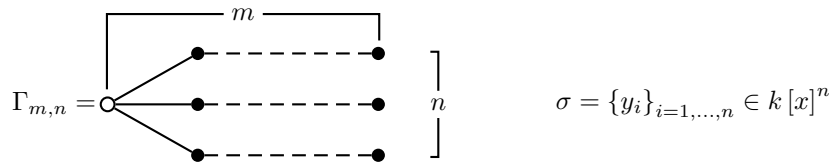
\begin{figure}[h]

\begin{pspicture}(-5,1.5)(10,-1)

\def\dedge{\ncline[linestyle=dashed]}

\rput(-2.5,0){$\Gamma_{m,n}=$}

\rput(0,0){

\pstree[treemode=R,radius=2.5pt,treesep=0.5cm,levelsep=1.2cm]{\Tc{3pt}}{

     \pstree{\TC*}{\skiplevels{1} \TC*[edge=\dedge] \endskiplevels}

     \pstree{\TC*}{\skiplevels{1} \TC*[edge=\dedge] \endskiplevels}

    \pstree{\TC*}{\skiplevels{1}  \TC*[edge=\dedge] \endskiplevels}

}

}

\pnode(-1.8,0.2){C}\pnode(1.82,0.85){D}

\ncbar[angleA=90]{C}{D}\ncput*[npos=1.5]{$m$}

\pnode(2.5,0.7){E}\pnode(2.5,-0.7){F}

\ncbar[arm=3pt]{E}{F}\ncput*[npos=1.5]{$n$}

\rput(6,0){$\sigma=\left\{y_i\right\}_{i=1,\ldots ,n}\in k\left[x\right]^n  $}

\end{pspicture}

\end{figure}

\noindent In Theorem (\ref{thm:Equivalence_of_Cats}), we establish
that the category $\overrightarrow{\mathfrak{D}}_{\left(A,x\right)}$
of DFS's with base $\left(A,x\right)$ equipped with certain birational
morphisms as above and the category $\mathcal{T}_{\left(A,x\right)}$
of $\left(A,x\right)$-labelled trees are equivalent. Then, in Theorem
(\ref{EquivalentTreesTheorem}), we classify these DFS's up to $X$-isomorphisms
in terms of the corresponding trees. In section $4$, we decompose
birational $X$-morphisms between two DFS's into a succession of simple
affine modifications \cite{KaZa99}, which we call \emph{fibered modifications}
(Theorem \ref{FiberedModifThm}). This leads to a canonical procedure
for constructing embeddings of a DFS $\pi:S\rightarrow X$ into a
projective $X$-scheme $\bar{\pi}:\bar{S}\rightarrow X$ (Proposition
\ref{pro:CanonicalCompletion}). Section $5$ is devoted to the study
of DFS's with base $\left(k\left[x\right],x\right)$, which we call
simply \emph{Danielewski surfaces}. We recall that the Makar-Limanov
invariant \cite{KML97} of an affine variety $V/k$ is defined as
the intersection in $\Gamma\left(V,\mathcal{O}_{V}\right)$ of all
the invariant rings of $\mathbb{G}_{a,k}$-actions on $V$. Makar-Limanov
\cite{ML90} and \cite{ML01} established that a surface $S_{P,m}\simeq S\left(\gamma_{P,m}\right)$
as above has trivial Makar-Limanov invariant $\textrm{ML}\left(S_{P,m}\right)=k$
if and only if $m=1$ and $\deg\left(P\right)\geq1$. More generally,
in case that $k=\bar{k}$ is algebraically closed, we give the following
characterization (Theorem \ref{MLGDSThm}).

\begin{thm}
A Danielewski surface $S\left(\gamma\right)$ has a trivial Makar-Limanov
invariant if and only if $\gamma$ is a comb, \emph{i.e} a labelled
tree $\gamma=\left(\Gamma,\sigma\right)$ with the property that all
but at most one of the direct descendants of a given $e\in\Gamma$
are terminal elements of $\Gamma$.
\end{thm}
\noindent Finally, we obtain the following description of normal
affine surfaces with a trivial Makar-Limanov invariant, which generalizes
previous results obtained by Daigle-Russell \cite{DaRus02} and Miyanishi-Masuda
\cite{MiMa03} for the particular case of log $\mathbb{Q}$-homology
planes (Theorem \ref{thm:MLCarac}). 

\begin{thm}
Every normal affine surface $S/\bar{k}$ with a trivial Makar-Limanov
invariant is isomorphic to a cyclic quotient of a Danielewski surface. 
\end{thm}

\section{Preliminaries on rooted trees}

\noindent Here we give some results on rooted trees and labelled
rooted trees which will be used in the following sections.

\subsection{Basic facts on rooted trees}

Let $\left(G,\leq\right)$ be a nonempty, finite, partially ordered
set (a \emph{poset}, for short). A totally ordered subset $C\subset G$
is called a \emph{subchain of length} $l\left(C\right)=\textrm{Card}\left(C\right)-1$.
\emph{A} subchain of maximal length is called a \emph{maximal subchain.}
For every $e\in G$, we let $\left(\uparrow e\right)_{G}=\left\{ e'\in G,e'\geq e\right\} $
and $\left(\downarrow e\right)_{G}=\left\{ e'\in G,e'\leq e\right\} $.
The \emph{edges} $E\left(G\right)$ of $G$ are the subsets $\overrightarrow{ee'}=\left\{ e'<e\right\} $
of $G$ with two elements such that $\left(\uparrow e'\right)_{G}\cap\left(\downarrow e\right)_{G}=\overrightarrow{ee'}$. 

\begin{defn}
A (\emph{rooted}) \emph{tree $\Gamma$} is poset with a unique minimal
element $e_{0}$ called the \emph{root}, and \emph{}such that $\left(\downarrow e\right)_{\Gamma}$
is a chain for every $e\in\Gamma$. A (\emph{rooted}) \emph{subtree}
of $\Gamma$ is a sub-poset $\Gamma'\subset\Gamma$ which is a rooted
tree for the induced order. We say that $\Gamma'$ is \emph{maximal}
if there exists $e'\in\Gamma$ such that $\Gamma'=\left(\uparrow e'\right)_{\Gamma}$. 
\end{defn}
\begin{enavant} An element $e\in\Gamma$ such that $l\left(\downarrow e\right)_{\Gamma}=m$
is said to be \emph{at level} $m$. The maximal elements $e_{i}=e_{i,m_{i}}$,
where $m_{i}=l\left(\downarrow e_{i}\right)_{\Gamma}$, of $\Gamma$
are called the \emph{leaves} of $\Gamma$. We denote the set of those
elements by $L\left(\Gamma\right)$. The corresponding subchains \begin{equation}
\left(\downarrow e_{i,m_{i}}\right)_{\Gamma}=\left\{ e_{i,0}=e_{0}<e_{i,1}<\cdots<e_{i,m_{i}-1}<e_{i,m_{i}}=e_{i}\right\} ,\quad i=1,\ldots,n.\label{eq:MaxChains}\end{equation}
 are the maximal subchains of $\Gamma$. We say that $\Gamma$ has
\emph{height} $h\left(\Gamma\right)=\max\left(m_{i}\right)$. The
\emph{first common ancestor} of two element $e,e'\in\Gamma\setminus\left\{ e_{0}\right\} $
is the maximal element of the chain $\left(\left(\downarrow e\right)_{\Gamma}\setminus\left\{ e\right\} \right)\cap\left(\left(\downarrow e'\right)_{\Gamma}\setminus\left\{ e'\right\} \right)$.
The \emph{children} of $e\in\Gamma\setminus L\left(\Gamma\right)$
are the minimal elements of $\left(\uparrow e\right)_{\Gamma}\setminus\left\{ e\right\} $.
We denote the set of those elements by $\textrm{Ch}_{\Gamma}\left(e\right)$. 

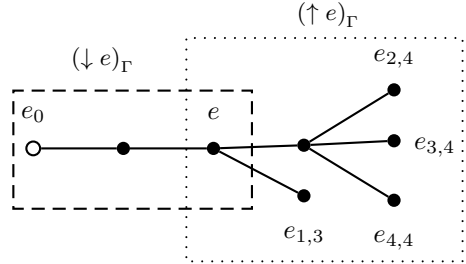
\begin{figure}[h]

\begin{pspicture}(-2,1.8)(10,-1.3)

\rput(4,0){

\pstree[treemode=R,radius=2.5pt,treesep=0.5cm,levelsep=1.2cm]{\Tc{3pt}~[tnpos=a]{$e_0$}}{

 \pstree{\TC*}{  \pstree{\TC*~[tnpos=a]{$e$}}   {\Tn \pstree{\TC*} {\TC*~[tnpos=a]{$e_{2,4}$}\TC*~[tnpos=r]{$e_{3,4}$}\TC*~[tnpos=b]{$e_{4,4}$}} \TC*~[tnpos=b]{$e_{1,3}$} }   }

}

}

\psframe[linestyle=dotted](3.3,-1.5)(7,1.5)

\rput(5.2,1.8){{\small $\left(\uparrow e \right)_\Gamma $}}

\psframe[linestyle=dashed](1,-0.8)(4.2,0.8)

\rput(2.2,1.2){{\small $\left(\downarrow e \right)_\Gamma $}}

\end{pspicture}

\caption{A tree  $\Gamma$ rooted in  $e_0$.}

\end{figure}

\end{enavant}

\begin{defn}
\label{BasicTreeMorphism} A \emph{morphism of} (\emph{rooted}) \emph{trees}
is an order-preserving map $\tau:\Gamma'\rightarrow\Gamma$ satisfying
the following properties: 

\noindent$\;$a) The image of a maximal chain of $\Gamma'$ by $\tau$
is a maximal chain of $\Gamma$.

\noindent$\;$b) For every $e'\in\Gamma'$, $\tau^{-1}\left(\tau\left(e'\right)\right)$
is either $e'$ itself or a maximal subtree of $\Gamma'$.

\noindent Injective and bijective morphisms are referred to as embeddings
and isomorphisms respectively. 
\end{defn}
\begin{enavant} \label{txt:Tree_Mor_Desc} \label{TreeIsoRemark}
This definition implies that for every leaf $e'_{i,m_{i}'}$ of  $\Gamma'$
at level $m'_{i}$, $\tau\left(e'_{i,m_{i}'}\right)$ is a leaf $e_{j\left(i\right),m_{j\left(i\right)}}$
of $\Gamma$ at level $m_{j\left(i\right)}\leq m'_{i}$, such that
$\tau\left(e'_{i,k}\right)=e_{j\left(i\right),\min\left(k,m_{j\left(i\right)}\right)}$
for every $k=0,\ldots,m'_{i}$ (see (\ref{eq:MaxChains}) for the
notation). 

\end{enavant}

\subsection{Labelled rooted trees}

In this subsection, we fix a pair $\left(A,x\right)$ consisting of
a domain $A$ and an element $x\in A\setminus\left\{ 0\right\} $. 

\begin{defn}
\label{def:W_Rooted_Tree} \label{def:W_Rooted_Tree_Mor} An $\left(A,x\right)$\emph{-labelled
tree} is \emph{}a pair $\gamma=\left(\Gamma,\sigma\right)$ consisting
of a a tree $\Gamma$ with the leaves $e_{1,m_{1}},\ldots,e_{n,m_{n}}$
and a cochain $\sigma=\left\{ \sigma_{i}\right\} _{i=1,\ldots,n}\in A^{n}$
such that $\sigma_{j}-\sigma_{i}\in x^{d_{ij}}A\setminus x^{d_{ij}+1}A$
whenever the first common ancestor of $e_{i,m_{i}}$ and $e_{j,m_{j}}$
is at level $d_{ij}<\min\left(m_{i},m_{j}\right)$. A cochain $\sigma$
with this this property is said to be $\Gamma$-\emph{compatible}.

A \emph{morphism of} $\left(A,x\right)$\emph{-labelled trees} $\tau:\left(\Gamma',\sigma'\right)\rightarrow\left(\Gamma,\sigma\right)$
is a morphism of trees $\tau:\Gamma'\rightarrow\Gamma$ such that
if $e_{j\left(i\right),m_{j\left(i\right)}}\in L\left(\Gamma\right)$
is the image of a leaf $e'_{i,m'_{i}}\in L\left(\Gamma'\right)$ by
$\tau$, then $\sigma_{i}-\sigma_{j\left(i\right)}\in x^{m_{j\left(i\right)}}A$.

The category of $\left(A,x\right)$-labelled just defined is denoted
by $\mathcal{T}_{\left(A,x\right)}$. 
\end{defn}
\begin{example}
\label{exa:Weighted_Rooted_Trees} Every isomorphism class of $\left(k\left[x\right],x\right)$-labelled
trees contains a tree $\gamma=\left(\Gamma,\sigma\right)$ with the
leaves $e_{1,m_{1}},\ldots,e_{n,m_{n}}$, such that $\sigma_{i}=\sum_{j=0}^{m_{i}-1}w_{i,j}x^{j}\in k\left[x\right]$
is a polynomial of degree $<m_{i}$ for every $i=1,\ldots,n$. In
turn, this cochain $\sigma$ is uniquely determined by the choice
of a \emph{weight function} \[
w_{\sigma}:E\left(\Gamma\right)\longrightarrow k,\quad\overrightarrow{e_{i,j+1}e_{i,j}}\mapsto w_{\sigma}\left(\overrightarrow{e_{i,j+1}e_{i,j}}\right)=w_{i,j}\]
 with the property that $w_{\sigma}\left(\overrightarrow{e'e}\right)\neq w_{\sigma}\left(\overrightarrow{e''e}\right)$
whenever $e'$ and $e''$ are children of the same $e\in\Gamma$.
A tree $\Gamma$ equipped with such a function $w$ is referred to
as a \emph{fine} $k$-\emph{weighted tree}. A morphism of trees $\tau:\Gamma'\rightarrow\Gamma$
gives rise to a morphism of $\left(k\left[x\right],x\right)$-labelled
$\tau:\left(\Gamma',\sigma'\right)\rightarrow\left(\Gamma,\sigma\right)$
if and only if $w_{\sigma}\left(\overrightarrow{\tau\left(e\right)\tau\left(e'\right)}\right)=w_{\sigma'}\left(\overrightarrow{ee'}\right)$whenever
$\tau\left(e\right)\neq\tau\left(e'\right)$. If it is the case, then
we say that $\tau$ is a \emph{morphism of} \emph{fine} $k$-\emph{weighted
trees}. In this way, we obtain a bijection between morphisms of $\left(A,x\right)$-labelled
trees and morphisms of fine $k$-weighted trees. We conclude that
the categories $\mathcal{T}_{\left(k\left[x\right],x\right)}$ of
$\left(k\left[x\right],x\right)$-labelled trees and $\mathcal{T}_{w}^{k}$
of fine $k$-weighted trees are equivalent. 
\end{example}

\subsection{Gluing trees}

Intuitively, $\left(A,x\right)$-labelled trees are constructed by
gluing $\left(A,x\right)$-labelled chains. Similarly, a morphism
$\tau:\gamma'\rightarrow\gamma$ of $\left(A,x\right)$-labelled trees
is uniquely determined by its values on maximal subchains of $\gamma'$.
More precisely, we have the following results. 

\begin{prop}
\label{Pro:TreeConstructProp} Suppose that $\left(n,m,d,\sigma\right)$
is a data consisting of an integer $n\geq1$, a multi-index $m=\left(m_{1},\ldots,m_{n}\right)\in\mathbb{Z}_{>0}^{n}$,
a matrix $d=\left(d_{ij}\right)_{i,j=1,\ldots,n}\in\textrm{Mat}_{n}\left(\mathbb{Z}_{\geq0}\right)$
and a cochain $\sigma=\left\{ \sigma_{i}\right\} _{i=1,\ldots,n}\in A^{n}$
with the following properties. 

1) For every $i\neq j$ , $d_{ij}=d_{ji}<\min\left(m_{i},m_{j}\right)$
. 

2) For every triple $i,j,k$, $\min\left(d_{ij},d_{ik}\right)=\min\left(d_{ji},d_{jk}\right)$. 

3) For every $i\neq j$, $\sigma_{j}-\sigma_{i}\in x^{d_{ij}}A\setminus x^{d_{ij}+1}A$. 

\noindent Then there exists a tree $\Gamma$, unique up to a unique
isomorphism, with the root $e_{0}$ and the leaves $e_{i}$ at levels
$m_{i}$, $i=1,\ldots,n$, such that $\sigma$ is $\Gamma$-compatible.
\end{prop}
\begin{proof}
Up to an isomorphism, $\gamma_{i}=\left(C_{i}=\left\{ e_{i,0}<e_{i,1}<\cdots<e_{i,m_{i}}\right\} ,\sigma_{i}\right)$
is the unique chain of height $m_{i}\geq1$ such that $\sigma_{i}$
is $C_{i}$-compatible. For every $i\neq j$, we let $C_{ij}=\left(\downarrow e_{i,d_{ij}}\right)_{C_{i}}\subset C_{i}$.
Condition (2) guarantees that there exist isomorphisms of chains $\phi_{ij}:C_{ij}\stackrel{\sim}{\rightarrow}C_{ij}$
such that $\phi_{ij}\left(C_{ij}\cap C_{ik}\right)=C_{ji}\cap C_{jk}$
for every triple $i,j,k$, and such that the cocycle condition $\phi_{ik}=\phi_{jk}\circ\phi_{ij}$
holds on on $C_{ij}\cap C_{ik}$. Moreover, (1) implies that for every
$i=1,\ldots,n$, $e_{i,0}$ belongs to $C_{ij}$ whereas $e_{i,m_{i}}$
does not. Therefore, the quotient poset $\Gamma=\bigsqcup_{i=1}^{n}C_{i}/\left(C_{ij}\ni e_{i,k}\sim\phi_{ij}\left(e_{j,k}\right)\in C_{ji}\right)$
is a tree rooted in the common image $e_{0}=\overline{e_{i,0}}$ of
the roots $e_{i,0}$ of the chains $C_{i}$, and with the leaves $e_{i}=\overline{e_{i,m_{i}}}$
at level $m_{i}$, $i=1,\ldots,n$ . Finally, (3) means exactly that
$\sigma$ is $\Gamma$-compatible. This completes the proof. 
\end{proof}
\begin{prop}
\label{Pro:MorphismConstructionProp} Let $\gamma=\left(\Gamma,\sigma\right)$
and $\gamma'=\left(\Gamma',\sigma'\right)$ be two $\left(A,x\right)$-labelled
trees. Suppose we are given a collection of morphisms of $\left(A,x\right)$-labelled
trees \[
\tau_{i}:\gamma'_{i}=\left(\Gamma'_{i},\sigma'_{i}\right)=\left(\left(\downarrow e'_{i}\right)_{\Gamma'},\sigma'_{i}\right)\longrightarrow\gamma,\quad\textrm{where }L\left(\Gamma'\right)=\left\{ e'_{1},\ldots,e'_{n'}\right\} ,\]
 restricting to a same morphism of trees $\tau_{ij}=\tau_{ji}:\Gamma'_{i}\cap\Gamma'_{j}\rightarrow\tau_{i}\left(\Gamma'_{i}\right)\cap\tau_{j}\left(\Gamma'_{j}\right)$
for every $i\neq j$. Then there exists a unique morphism of $\left(A,x\right)$-labelled
trees $\tau:\gamma'\rightarrow\gamma$ such that $\tau\mid_{\gamma'_{i}}=\tau_{i}$
for every $i=1,\ldots,n'$.
\end{prop}
\begin{proof}
The conditions guarantee that there exists a unique order-preserving
map $\tau:\Gamma'\rightarrow\Gamma$ such that $\tau\mid_{\Gamma_{i}'}=\tau_{i}$,
and such that the preimage $\tau^{-1}\left(e\right)$ of a given $e\in\tau\left(\Gamma'\right)$
is the maximal subtree of $\Gamma'$ rooted in the unique minimal
element of $\tau^{-1}\left(e\right)$. Thus $\tau$ is a morphism
of trees compatible with the cochains $\sigma'$ and $\sigma$ in
the sense (\ref{def:W_Rooted_Tree_Mor}), whence a morphism of $\left(A,x\right)$-labelled
trees.
\end{proof}

\subsection{Blow-downs of trees}

By definition, (see (\ref{def:W_Rooted_Tree_Mor})), A morphism of
trees $\tau:\Gamma'\rightarrow\Gamma$ factors through the retraction
of a collection of maximal subtrees of $\Gamma'$, followed by an
embedding. Therefore, for every element $e\in\Gamma'\setminus L\left(\Gamma'\right)$
such that $\textrm{Ch}_{\Gamma'}\left(e\right)\subset L\left(\Gamma'\right)$,
the image of the subtree $\Gamma'_{e}=\left(\uparrow e\right)_{\Gamma'}$
of $\Gamma'$ is either a subtree $\tau\left(\Gamma'_{e}\right)$
of $\Gamma$ isomorphic to $\Gamma'_{e}$ or the unique element $\tau\left(e\right)\in L\left(\Gamma\right)$.
In the second case, $\tau$ factors through the morphism of trees
\[
\tau_{e}:\Gamma'\longrightarrow\Gamma'\setminus\textrm{Ch}_{\Gamma'}\left(e\right)\quad e'\mapsto\left\{ \begin{array}{rl}
e' & \textrm{if }e'\in\Gamma'\setminus\textrm{Ch}_{\Gamma'}\left(e\right)\\
e & \textrm{if }e'\in\textrm{Ch}_{\Gamma'}\left(e\right)\end{array}\right.\]

\begin{defn}
\label{def:TreeBlowDef} Let $\gamma=\left(\Gamma,\sigma\right)$
be an $\left(A,x\right)$-labelled tree, and let $e\in\Gamma\setminus L\left(\Gamma\right)$
be an element such that $\textrm{Ch}_{\Gamma}\left(e\right)\subset L\left(\Gamma\right)$.
A \emph{blow-down of the leaves at} $e$ is an $\left(A,x\right)$-labelled
tree $\hat{\gamma}\left(e\right)=\left(\hat{\Gamma}\left(e\right),\hat{\sigma}\left(e\right)\right)$
with underlying tree $\hat{\Gamma}\left(e\right)=\Gamma\setminus\textrm{Ch}_{\Gamma}\left(e\right)$
for which the morphism of trees $\tau_{e}$ above is a morphism of
$\left(A,x\right)$-labelled trees $\tau_{e}:\gamma\rightarrow\hat{\gamma}\left(e\right)$.
Since two labelled trees $\hat{\gamma}\left(e\right)$ with this property
are isomorphic, the morphism $\tau_{e}$ itself will be usually referred
to as \emph{the} blow-down of the leaves at $e$. 
\end{defn}
\noindent As a consequence of the above discussion, we obtain the
following description. 

\begin{prop}
\label{Pro:TreeBlowUpLemma} A morphism of $\left(A,x\right)$-labelled
trees factors into a sequence of blow-downs of leaves followed by
an embedding.
\end{prop}

\subsection{Equivalence of labelled trees}

An $\left(A,x\right)$-labelled tree $\gamma=\left(\Gamma,\sigma\right)$
rooted in $e_{0}$ is called \emph{essential} if $\textrm{Card}\left(\textrm{Child}_{\Gamma}\left(e_{0}\right)\right)\neq1$.
The \emph{essential subtree} $\textrm{Es}\left(\Gamma\right)$ of
a given tree $\Gamma$ is the maximal subtree of $\Gamma$ rooted
either in the unique leaf of $\Gamma$ if $\Gamma$ is a chain, or
in the first common ancestor $\tilde{e}_{0}$ of the leaves of $\Gamma$
otherwise. For instance, the tree $\Gamma$ of Figure 1 above is not
essential, and its essential rooted subtree is the maximal subtree
$\left(\uparrow e\right)_{\Gamma}$ of $\Gamma$ rooted in $e$. If
an $\left(A,x\right)$-labelled tree $\gamma=\left(\Gamma,\sigma\right)$
is not essential then there exists $c\in A$ and an $\textrm{Es}\left(\Gamma\right)$-compatible
cochain $\textrm{Es}\left(\sigma\right)=\left\{ \tilde{\sigma}_{i}\right\} _{i=1,\ldots,n}\in A^{n}$
such that $\sigma_{i}=c+x^{m}\tilde{\sigma}_{i}$ for every $i=1,\ldots,n$,
where $m$ denotes the height of the root $\tilde{e}_{0}$ of the
essential subtree $\textrm{Es}\left(\Gamma\right)$ of $\Gamma$.
A cochain $\textrm{Es}\left(\sigma\right)$ with this property is
called an \emph{essential cochain} for $\gamma$, and we say that
$\textrm{Es}\left(\gamma\right)=\left(\textrm{Es}\left(\Gamma\right),\textrm{Es}\left(\sigma\right)\right)$
is an \emph{essential labelled subtree for} $\gamma$.

\begin{defn}
\label{def:Equivalence_of_trees} We say that two $\left(A,x\right)$-labelled
trees $\gamma=\left(\Gamma,\sigma\right)$ and $\gamma'=\left(\Gamma',\sigma'\right)$
are \emph{equivalent} if there exist essential cochains $\textrm{Es}\left(\sigma\right)$
and $\textrm{Es}\left(\sigma'\right)$ for $\gamma$ and $\gamma'$
respectively, an isomorphism of trees $\tau:\textrm{Es}\left(\Gamma'\right)\stackrel{\sim}{\rightarrow}\textrm{Es}\left(\Gamma\right)$
and a pair $\left(a,b\right)\in A^{*}\times A$ such that \[
a\textrm{Es}\left(\sigma'\right)_{i}-\textrm{Es}\left(\sigma\right)_{j\left(i\right)}+b\in x^{m_{j\left(i\right)}}A\]
 whenever $e_{j\left(i\right),m_{j\left(i\right)}}\in L\left(\textrm{Es}\left(\Gamma\right)\right)$
is the image of $e'_{i,m_{i}'}\in L\left(\textrm{Es}\left(\Gamma'\right)\right)$
by $\tau$.
\end{defn}
\begin{example}
By definition, two essential cochains for a given labelled tree $\gamma$
differ by the addition of a constant $b\in A$. Therefore, the essential
labelled subtrees for $\gamma$ are all equivalent in the sense of
definition (\ref{def:Equivalence_of_trees}). 
\end{example}

\section{Danielewski-Fieseler surfaces as fiber bundles }

To fix the notation, we let $k$ be a field of caracteristic zero,
and we let $\left(A,x\right)$ be a pair consisting either of a discrete
valuation ring $A$ with uniformizing parameter $x$ and residue field
$k$, or a polynomial ring $k\left[x\right]$ in one variable $x$.
We let $X=\textrm{Spec}\left(A\right)$, $x_{0}=\textrm{div}\left(x\right)\in X$
, and we denote by $X_{*}=X\setminus\left\{ x_{0}\right\} \simeq\textrm{Spec}\left(A_{x}\right)$
the open complement of $x_{0}$ in $X$. 

\begin{defn}
\label{def:DFS_def} A \emph{Danielewski-Fieseler surface} (a DFS
for short) with base $\left(X,x_{0}\right)$ (or, equivalently, with
base $\left(A,x\right)$) is an integral affine $X$-scheme $\pi:S\rightarrow X$
of finite type such that $S\mid_{X_{*}}$ is isomorphic to the trivial
line bundle $\ba_{X_{*}}^{1}$ over $X_{*}$ and such that the scheme-theoretic
fiber $\pi^{-1}\left(x_{0}\right)$ is nonempty and reduced, consisting
of a disjoint union of curves isomorphic to the affine line $\ba_{k}^{1}$.
A DFS with base $\left(k\left[x\right],x\right)$ is simply referred
to as a \emph{Danielewski surface}.

A \emph{morphism of Danielewski-Fieseler surfaces} with base $\left(X,x_{0}\right)$
is a birational $X$-morphism $\beta:S'\rightarrow S$ restricting
to an $X_{*}$-isomorphism $\beta_{*}:S'\mid_{X_{*}}\stackrel{\sim}{\longrightarrow}S\mid_{X_{*}}$. 

Danielewski-Fieseler surfaces together with these morphisms form a
sub-category $\mathfrak{D}_{\left(A,x\right)}$ (or $\mathfrak{D}_{\left(X,x_{0}\right)}$)
of the category $\left(\textrm{Sch}_{/X}\right)$ of $X$-schemes. 
\end{defn}
\begin{enavant} The total space $S$ of DFS $\pi:S\rightarrow X$
is a smooth affine surface over $k$. Indeed, $X/k$ is itself affine
and smooth, and the local criterion for flatness (\cite[III,§5]{Bour})
guarantees that $\pi$ is a faithfully flat morphism, whence a smooth
morphism as its geometric fibers are regular. Note that in contrast
with general DFS's, a Danielewski surface is an affine $k$-scheme
of \emph{finite type} since in this case $X\simeq\ba_{k}^{1}$ by
definition. In general, an arbitrary $\ba^{1}$-fibration $\pi:S\rightarrow\ba_{k}^{1}$
on a smooth affine surface $S$ does not give rise to a structure
of Danielewski surface on $S$, but this does not prevent $S$ from
being a Danielewski surface for another suitable $\ba^{1}$-fibration.
For instance, the Bandman and Makar-Limanov surface \cite{BML01}
$S\subset\textrm{Spec}\left(k\left[x\right]\left[y,z,u\right]\right)$
with equations \[
xz-y\left(y-1\right)=0,\quad yu-z^{2}=0,\quad xu-\left(y-1\right)z=0,\]
 equipped with the $\ba^{1}$-fibration $pr_{u}:S\rightarrow\textrm{Spec}\left(k\left[u\right]\right)$
is not a Danielewski surface with base $\left(k\left[u\right],u\right)$.
Indeed the fiber $pr_{u}^{-1}\left(0\right)$ is not reduced. However,
it is a Danielewski surface $pr_{x}:S\rightarrow\textrm{Spec}\left(k\left[x\right]\right)$
with base $\left(k\left[x\right],x\right)$. 

\end{enavant} 

\noindent Following Fieseler \cite{Fie94}, DFS's with base $\left(X,x_{0}\right)$
can be described as certain $\ba^{1}$-bundles $\rho:S\rightarrow Y$
over a curve $\delta:Y\rightarrow X$ with an $n$-fold point over
$x_{0}$. In this section, we give a new interpretation of Fieseler's
'cocycle construction'.

\subsection{Torsors under a line bundle}

Given an invertible sheaf $\mathcal{L}$ on a scheme $X$, the line
bundle \[
\mathbf{V}\left(\mathcal{L}\right)=\mathbf{Spec}\left(\mathbf{S}\left(\mathcal{L}\right)\right)=\mathbf{Spec}\left(\bigoplus_{n\geq0}\mathcal{L}^{n}\right)\longrightarrow X\]
 is equipped with the structure of a commutative group scheme for
the group law $m:\mathbf{V}\left(\mathcal{L}\right)\times_{X}\mathbf{V}\left(\mathcal{L}\right)\rightarrow\mathbf{V}\left(\mathcal{L}\right)$
defined by means of the $\mathcal{O}_{X}$-algebras homomorphism $\mathbf{S}\left(\mathcal{L}\right)\longrightarrow\mathbf{S}\left(\mathcal{L}\right)\otimes\mathbf{S}\left(\mathcal{L}\right)$
induced by the diagonal homomorphism $\Delta:\mathcal{L}\rightarrow\mathcal{L}\oplus\mathcal{L}$.
For instance, if $\mathcal{L}=\mathcal{O}_{X}$ then $\mathbf{V}\left(\mathcal{L}\right)$
is simply the additive group scheme $\mathbb{G}_{a,X}=\mathbf{Spec}\left(\mathcal{O}_{X}\left[T\right]\right)$. 

\begin{enavant} A scheme $\pi:S\rightarrow X$ locally isomorphic
to the trivial line bundle $\ba_{X}^{1}$ is called an $\ba^{1}$\emph{-bundle}.
If $S$ comes further equipped with an action of a line bundle $\mathbf{V}\left(\mathcal{L}\right)/X$
for which there exists an open covering $\mathcal{U}=\left(X_{i}\right)_{i\in I}$
of $X$ such that $S\mid_{X_{i}}$ is equivariantly isomorphic to
$\mathbf{V}\left(\mathcal{L}\right)\mid_{X_{i}}$ acting on itself
by translations for every $i\in I$, then we say that $S$ is a $\mathbf{V}\left(\mathcal{L}\right)$-\emph{torsor}.
We recall \cite[XI.4.7]{SGA1} that the set of isomorphism classes
of $\mathbf{V}\left(\mathcal{L}\right)$-torsors is a group isomorphic
to $H^{1}\left(X,\mathcal{L}^{\vee}\right)$. Indeed, an automorphism
of the symmetric $\mathcal{O}_{X}$-algebra $\mathbf{S}\left(\mathcal{L}\right)$
is equivariant for the co-action of $\mathbf{S}\left(\mathcal{L}\right)$
on itself by translations if and only if it is induced by an homomorphism
of $\mathcal{O}_{X}$-modules $\left(g,\textrm{Id}\right):\mathcal{L}\rightarrow\mathcal{O}_{X}\oplus\mathcal{L}$,
where $g\in\textrm{Hom}_{X}\left(\mathcal{L},\mathcal{O}_{X}\right)$.
Therefore, given a $\mathbf{V}\left(\mathcal{L}\right)$-torsor $\pi:S\rightarrow X$
which becomes trivial on $\mathcal{U}$, there exists a \v{C}ech
cocycle $g=\left\{ g_{ij}\right\} _{i,j\in I}\in C^{1}\left(\mathcal{U},\mathcal{L}^{\vee}\right)$
such that $S$ is isomorphic to the scheme $\mathbf{W}\left(\mathcal{U},\mathcal{L},g\right)$
defined as the spectrum of the quasi-coherent $\mathcal{O}_{X}$-algebra
obtained by gluing the symmetric algebras $\mathbf{S}\left(\mathcal{L}\mid_{X_{i}}\right)$
over $X_{ij}=X_{i}\cap X_{j}$ by means of the $\mathcal{O}_{X_{ij}}$-algebras
isomorphisms $\mathbf{S}\left(g_{ij},\textrm{Id}\right)$ induced
by the $\mathcal{O}_{X_{ij}}$-modules homomorphisms $\left(g_{ij},\textrm{Id}\right)\in\textrm{Hom}_{X_{ij}}\left(\mathcal{L}\mid_{X_{ij}},\mathcal{O}_{X_{ij}}\oplus\mathcal{L}\mid_{X_{ij}}\right)$,
$i\neq j$. The isomorphism class of $\mathbf{W}\left(\mathcal{U},\mathcal{L},g\right)$
is simply the image in $H^{1}\left(X,\mathcal{L}^{\vee}\right)\simeq\check{H}^{1}\left(X,\mathcal{L}^{\vee}\right)$
of the class $\left[g\right]\in\check{H}^{1}\left(\mathcal{U},\mathcal{L}^{\vee}\right)$
of $g\in C^{1}\left(\mathcal{U},\mathcal{L}^{\vee}\right)$. 

\end{enavant}

\begin{enavant} \label{pro:Torsor_equiv_carac} If the base scheme
$X$ is integral, then every $\ba^{1}$-bundle $\pi:S\rightarrow X$
admits a structure of a $\mathbf{V}\left(\mathcal{L}\right)$-torsor
for a suitable invertible sheaf $\mathcal{L}$ on $X$. Indeed, since
$X$ is integral, the transition isomorphisms $\tau_{ij}=\tau_{i}\circ\tau_{j}^{-1}$
associated with a given collection of trivialisations $\tau_{i}:S\mid_{X_{i}}\stackrel{\sim}{\rightarrow}\textrm{Spec}\left(\mathcal{O}_{X_{i}}\left[T\right]\right)$,
$i\in I$, are induced by automorphisms $T\mapsto\tilde{g}_{ij}+f_{ij}T$
of $\mathcal{O}_{X_{ij}}\left[T\right]$, where $\tilde{g}_{ij}\in\Gamma\left(X_{ij},\mathcal{O}_{X}\right)$
and $f_{ij}\in\Gamma\left(X_{ij},\mathcal{O}_{X}^{*}\right)$. For
a triple of indices $i\neq j\neq k$, the identity $\tau_{ik}=\tau_{jk}\circ\tau_{ij}$
over $X_{i}\cap X_{j}\cap X_{k}$ guarantees that $\left(f_{ij}\right)_{i,j\in I}\in C^{1}\left(\mathcal{U},\mathcal{O}_{X}^{*}\right)$
is a \v{C}ech cocycle defining an invertible sheaf $\mathcal{L}$
on $X$, trivial on $\mathcal{U}$, with isomorphisms $\phi_{i}:\mathcal{L}\mid_{X_{i}}\stackrel{\sim}{\rightarrow}\mathcal{O}_{X_{i}}$,
$i\in I$. By construction, $g=\left\{ \tilde{g}_{ij}\cdot\phi_{i}\mid_{X_{ij}}\right\} _{i,j\in I}\in C^{1}\left(\mathcal{U},\mathcal{L}^{\vee}\right)$
is a \v{C}ech cocycle for which $S$ is $X$-isomorphic to the $\mathbf{V}\left(\mathcal{L}\right)$-torsor
$\mathbf{W}\left(\mathcal{U},\mathcal{L},g\right)$. 

\end{enavant}

\subsection{Schemes with an $n$-fold divisor}

Given a principal divisor $x_{0}=\textrm{div}\left(x\right)$ on an
integral scheme $X$ and an integer $n\geq1$, we let $\delta_{n}:X\left(n\right)=X\left(x_{0},n\right)\rightarrow X$
be the scheme obtained by gluing $n$ copies $d_{i}:X_{i}\left(n\right)\stackrel{\sim}{\rightarrow}X$
of $X$ by the identity over the open subsets $X_{i}\left(n\right)_{*}=X_{i}\left(n\right)\setminus x_{i}\left(n\right)$,
where $x_{i}\left(n\right)=d_{i}^{-1}\left(x_{0}\right)$. We denote
by $\mathcal{U}_{\left(n\right)}=\left(X_{i}\left(n\right)\right)_{i=1,\ldots,n}$
the canonical open covering of $X\left(n\right)$. For every $i\neq j$,
we let $X_{ij}\left(n\right)=X_{i}\left(n\right)\cap X_{j}\left(n\right)\simeq X\setminus\left\{ x_{0}\right\} $. 

\begin{defn}
\label{PicXPrelim} Given a multi-index $\mu=\left(\mu_{1},\ldots,\mu_{n}\right)\in\mathbb{Z}^{n}$,
we let \begin{eqnarray*}
\mathcal{L}_{n,\mu} & = & \mathcal{O}_{X\left(n\right)}\left(\mu_{1}x_{1}\left(n\right)+\ldots+\mu_{n}x_{n}\left(n\right)\right)\end{eqnarray*}
 be sub-$\mathcal{O}_{X\left(n\right)}$-module of the constant sheaf
$\mathcal{K}_{X\left(n\right)}$ of rational functions on $X\left(n\right)$
generated by $\left(x\circ\delta_{n}\right)^{-\mu_{i}}$ on $X_{i}\left(n\right)$.
The dual sheaf $\mathcal{L}_{n,\mu}^{\vee}$ of $\mathcal{L}_{n,\mu}$
is isomorphic to the sheaf $\mathcal{L}_{n,-\mu}$ corresponding to
the multi-index $-\mu=\left(-\mu_{1},\ldots,-\mu_{n}\right)\in\mathbb{Z}^{n}$.
We denote by $s_{\mu}$ the canonical \emph{rational} section of $\mathcal{L}_{n,\mu}$
corresponding to the constant section $1\in\Gamma\left(X\left(n\right),\mathcal{K}_{X\left(n\right)}\right)$. 
\end{defn}

\subsection{Danielewski-Fieseler surfaces as torsors}

Given a DFS $\pi:S\rightarrow X$ with base $\left(X,x_{0}\right)$,
we denote by $C_{1},\ldots,C_{n}$ the connected component of the
fiber $\pi^{-1}\left(x_{0}\right)$. We let $\delta_{n}:X\left(n\right)=X\left(x_{0},n\right)\rightarrow X$
be as above, and we let $\rho_{n}:S\rightarrow X\left(n\right)$ be
the unique $X$-morphism such that $\pi=\delta_{n}\circ\rho_{n}$,
and such that $C_{i}=\rho_{n}^{-1}\left(x_{i}\left(n\right)\right)$
for every $i=1,\ldots,n$. 

\begin{prop}
\label{pro:DFS_as_torsors} There exist a multi-index $\mu=\left(\mu_{1},\ldots,\mu_{n,}\right)\in\mathbb{Z}_{\geq0}^{n}$
such that $\mu_{i}=0$ for at least one indice $i$, and a cocycle
$g\in C^{1}\left(\mathcal{U}_{\left(n\right)},\mathcal{L}_{n,\mu}\right)$
such that $S$ is $X\left(n\right)$-isomorphic to the $\mathbf{V}\left(\mathcal{L}_{n,-\mu}\right)$-torsor
$\rho_{n}:\mathbf{W}\left(\mathcal{U}_{\left(n\right)},\mathcal{L}_{n,-\mu},g\right)\rightarrow X\left(n\right)$.
\end{prop}
\begin{proof}
Since $X$ is either the affine line or the spectrum of a discrete
valuation ring, the Picard group $\textrm{Pic}\left(X\right)$ is
trivial, and so every invertible sheaf on $X\left(n\right)$ is isomorphic
to $\mathcal{L}_{n,-\mu}$ for a certain multi-index $\mu\in\mathbb{Z}_{\geq0}^{n}$
such that $\mu_{i}=0$ for at least one indice $i$. Moreover, every
open subset $X_{i}\left(n\right)\simeq X$ of the covering $\mathcal{U}_{\left(n\right)}$
is affine. Therefore, a $\mathbf{V}\left(\mathcal{L}_{n,-\mu}\right)$-torsor
$\rho:W\rightarrow X\left(n\right)$ is isomorphic to $\mathbf{W}\left(\mathcal{U}_{\left(n\right)},\mathcal{L}_{n,-\mu},g\right)$
for a certain cocycle $g\in C^{1}\left(\mathcal{U},\mathcal{L}_{n,\mu}\right)$
representing its isomorphism class $c\left(W\right)\in H^{1}\left(X\left(n\right),\mathcal{L}_{n,\mu}\right)$.
So, by (\ref{pro:Torsor_equiv_carac}), it suffices to show that $\rho_{n}:S\rightarrow X\left(n\right)$
is an $\ba^{1}$-bundle . This can be done in a similar way as in
Lemma 1.2 in \cite{Fie94}. 
\end{proof}
\begin{enavant} \label{txt:DFS_carac_equiv} In view of the correspondence
(\ref{pro:Torsor_equiv_carac}) between $\ba^{1}$-bundles and $\mathbf{V}\left(\mathcal{L}\right)$-torsors,
a DFS $\mathbf{W}\left(\mathcal{U}_{\left(n\right)},\mathcal{L}_{n,-\mu},g\right)$
is $X$-isomorphic to the surface $S=\bigsqcup_{i=1}^{n}\ba_{X_{i}}^{1}/\sim$
obtained by gluing $n$ copies of $\ba_{X_{i}}^{1}=\textrm{Spec}\left(A\left[T_{i}\right]\right)$
of $\ba_{X}^{1}$ by means of $A_{x}$-algebras isomorphisms $A_{x}\left[T_{i}\right]\stackrel{\sim}{\rightarrow}A_{x}\left[T_{j}\right]$,
$T_{i}\mapsto\tilde{g}_{ij}+x^{\mu_{i}-\mu_{j}}T_{j}$, where $\tilde{g}_{ij}:=x^{\mu_{i}}g_{ij}\in A_{x}$
for every $i\neq j$. In this way, we recover the ''cocycle construction''
of Fieseler \cite{Fie94} . 

\end{enavant}

\begin{enavant}\label{txt:Affiness_criterion} A general scheme $S=\mathbf{W}\left(\mathcal{U}_{\left(n\right)},\mathcal{L}_{n,-\mu},g\right)$
is not affine. For instance, if $n\geq2$, then $\ba_{X\left(n\right)}^{1}$
is not even separated. However, a similar argument as in Proposition
1.4 in \cite{Fie94} shows that $S$ is affine if and only if the
corresponding \emph{transition functions} $\tilde{g}_{ij}=x^{\mu_{i}}g_{ij}$
and $\tilde{g}_{ji}=-x^{\mu_{j}}g_{ij}$ above have a pole in $x_{0}$
for every $i\neq j$. This means equivalently that for every $i\neq j$,
$g_{ij}\in\Gamma\left(X_{ij}\left(n\right),\mathcal{L}_{n,\mu}\right)$
is not in the image of the differential \[
\left(\textrm{res}_{i}-\textrm{res}_{j}\right):\Gamma\left(X_{j}\left(n\right),\mathcal{L}_{n,\mu}\right)\oplus\Gamma\left(X_{i}\left(n\right),\mathcal{L}_{n,\mu}\right)\rightarrow\Gamma\left(X_{ij}\left(n\right),\mathcal{L}_{n,\mu}\right).\]
 In turn, this condition is satisfied if and only if for every $i\neq j$,
$S\mid_{X_{i}\left(n\right)\cup X_{j}\left(n\right)}$ is a nontrivial
$\mathbf{V}\left(\mathcal{L}_{n,-\mu}\right)\mid_{X_{i}\left(n\right)\cup X_{j}\left(n\right)}$-torsor.
This leads to the following criterion. 

\end{enavant}

\begin{prop}
\label{pro:DFSCriterion} A scheme $\mathbf{W}\left(\mathcal{U}_{\left(n\right)},\mathcal{L}_{n,-\mu},g\right)$
is a DFS if and only if it restricts to a nontrivial $\mathbf{V}\left(\mathcal{L}_{n,-\mu}\right)\mid_{U}$-torsor
on every open subset $U\simeq X\left(2\right)$ of $X$. 
\end{prop}

\subsection{Morphisms of Danielewski-Fieseler surfaces}

Given two DFS's $S'=\mathbf{W}\left(\mathcal{U}_{\left(n'\right)},\mathcal{L}_{n',-\mu'},g'\right)$
and $S'=\mathbf{W}\left(\mathcal{U}_{\left(n'\right)},\mathcal{L}_{n',-\mu'},g'\right)$
as in (\ref{pro:DFS_as_torsors}) and a morphism of DFS's $\beta:S'\rightarrow S$,
we let $\alpha:X\left(n'\right)\rightarrow X\left(n\right)$ be the
unique $X$-morphism such that $\beta\left(\rho_{n'}^{-1}\left(x_{i}\left(n'\right)\right)\right)\subset\rho_{n}^{-1}\left(\alpha\left(x_{i}\left(n\right)\right)\right)$
for every $i=1,\ldots,n'$. We denote by $\bar{g}\in C^{1}\left(\mathcal{U}_{\left(n'\right)},\alpha^{*}\mathcal{L}_{n,\mu}\right)$
the image of $\alpha^{*}g\in C^{1}\left(\alpha^{-1}\left(\mathcal{U}_{\left(n\right)}\right),\alpha^{*}\mathcal{L}_{n,\mu}\right)$
by the restriction maps between \v{C}ech complexes. We obtain a factorization
$\beta=pr_{S}\circ\beta'$ \[\xymatrix{ S'\ar[rr]^-{\beta'} \ar[d]_{\rho_{n'}} & & \tilde{S}=S\times_{X\left(n\right)}X\left(n'\right)\simeq\mathbf{W}\left(\mathcal{U}_{\left(n'\right)},\alpha^*\mathcal{L}_{n,-\mu},\bar{g}\right) \ar[rr]^-{pr_S} \ar[d]_{p_2} && S \ar[d]_{\rho_n} \\ X\left(n'\right) \ar@{=}[rr] && X\left(n'\right) \ar[rr]^{\alpha} && X\left(n\right) }\]  
where $\beta'$ is an $X\left(n'\right)$-morphism restricting an
isomorphism over $\delta_{n'}^{-1}\left(x_{0}\right)$. This is the
case if and only if the homomorphisms $\mathbf{S}\left(\alpha^{*}\mathcal{L}_{n,-\mu}\right)\mid_{X_{i}\left(n'\right)}\rightarrow\mathbf{S}\left(\mathcal{L}_{n',-\mu'}\right)\mid_{X_{i}\left(n'\right)}$
corresponding to the restrictions $\beta'_{i}:S'\mid_{X_{i}\left(n'\right)}\rightarrow\tilde{S}\mid_{X_{i}\left(n'\right)}$
of $\beta'$ are induced by $\mathcal{O}_{X_{i}}$-modules homomorphisms
$\left(c_{i},\zeta_{i}\right):\alpha^{*}\mathcal{L}_{n,-\mu}\mid_{X_{i}\left(n'\right)}\rightarrow\mathcal{O}_{X_{i}\left(n'\right)}\oplus\mathcal{L}_{n',-\mu'}\mid_{X_{i}\left(n'\right)}$,
$i=1,\ldots,n'$, where $\zeta_{i}\in\textrm{Hom}_{X_{i}\left(n'\right)}\left(\alpha^{*}\mathcal{L}_{n,-\mu},\mathcal{L}_{n',-\mu'}\right)$
restricts to an isomorphism over $X_{i}\left(n'\right)\setminus\left\{ x_{i}\left(n'\right)\right\} $,
and where $c_{i}\in\Gamma\left(X_{i}\left(n'\right),\alpha^{*}\mathcal{L}_{n,\mu}\right)$.
Moreover, it follows from the definition of $S'$ and $\tilde{S}$
respectively that the $\beta_{i}$'s coincide on the overlaps $X_{ij}\left(n'\right)$,
$i\neq j$, if and only if the $\zeta_{i}$'s glue to a global section
$\zeta\in\textrm{Hom}_{X\left(n'\right)}\left(\alpha^{*}\mathcal{L}_{n,-\mu},\mathcal{L}_{n',-\mu'}\right)$
such that $g'\circ\zeta=\bar{g}+\partial\left(c\right)$, where $\partial$
denotes the differential of the \v{C}ech complex $C^{\bullet}\left(\mathcal{U}_{\left(n'\right)},\alpha^{*}\mathcal{L}_{n,-\mu}\right)$.
Summing up, we obtain the following characterization. 

\begin{prop}
\label{Pro:MorphismDescription} A morphism of DFS's $\beta:S'\rightarrow S$
exists if and only if there exists a data $\left(\alpha,\theta,c\right)$
consisting of: 

1) an $X$-morphism $\alpha:X\left(n'\right)\rightarrow X\left(n\right)$,

2) a section $\theta=^{t}\zeta\in\textrm{Hom}_{X\left(n'\right)}\left(\mathcal{L}_{n',\mu'},\alpha^{*}\mathcal{L}_{n,\mu}\right)$
such that $\textrm{Supp}\left(\theta\right)\subset\delta_{n'}^{-1}\left(x_{0}\right)$,

3) a cochain $c\in C^{0}\left(\mathcal{U}_{\left(n'\right)},\alpha^{*}\mathcal{L}_{n,\mu}\right)$
such that $\theta\left(g'\right)=\bar{g}+\partial\left(c\right)$. 

\noindent Furthermore, $\beta$ is an isomorphism if and only if
$\alpha$ and $\theta$ are. 
\end{prop}
\begin{example}
\label{exa:Mor_to_A1} Given a DFS $S=\mathbf{W}\left(\mathcal{U}_{\left(n\right)},\mathcal{L}_{n,-\mu},g\right)$,
there exists an integer $h_{0}\geq\max\left(\mu_{i}\right)$, depending
on $g$ and $\mu$, such that for every $h\geq h_{0}$, $x^{h}s_{-\mu}$
is a regular section of $\mathcal{L}_{n,-\mu}$ defining an homomorphism
of $\mathcal{O}_{X\left(n\right)}$-modules $\theta_{\mu,h}:\mathcal{L}_{n,\mu}\rightarrow\mathcal{O}_{X\left(n\right)}$
with the property that $\theta_{\mu,h}\left(g_{1j}\right)\in\Gamma\left(X_{1j}\left(n\right),\mathcal{O}_{X\left(n\right)}\right)$
extends to a section $\sigma_{j}$ of $\mathcal{O}_{X_{j}\left(n\right)}$
for every $j=2,\ldots,n$. By (\ref{Pro:MorphismDescription}), the
data $\left(\delta_{n},\theta_{\mu,h},\sigma\right)$, where $\sigma=\left\{ \sigma_{1}=0,\sigma_{j}\right\} _{j=2,\ldots,n}\in A^{n}$,
corresponds to a morphism of DFS's $\psi:S\rightarrow\ba_{X}^{1}$. 
\end{example}

\subsection{Additive group actions on Danielewski-Fieseler surfaces}

As a torsor, a DFS $\pi:S=\mathbf{W}\left(\mathcal{U}_{\left(n\right)},\mathcal{L}_{n,-\mu},g\right)\stackrel{\rho}{\longrightarrow}X\left(n\right)\stackrel{\delta_{n}}{\rightarrow}X$
comes equipped with an action $m_{n,\mu}:\mathbf{V}\left(\mathcal{L}_{n,-\mu}\right)\times_{X\left(n\right)}S\rightarrow S$
of the line bundle $\mathbf{V}\left(\mathcal{L}_{n,-\mu}\right)$.
Every nonzero section $s\in\Gamma\left(X\left(n\right),\mathcal{L}_{n,\mu}\right)$
gives rise to a morphism of group schemes $\phi_{s}:\mathbb{G}_{a,X\left(n\right)}\rightarrow\mathbf{V}\left(\mathcal{L}_{n,-\mu}\right)$,
whence to a nontrivial action \[
m_{s}=m_{n,\mu}\circ\left(\phi_{s}\times\textrm{Id}\right):\mathbb{G}_{a,X\left(n\right)}\times_{X\left(n\right)}S\stackrel{\phi_{s}\times\textrm{Id}}{\longrightarrow}\mathbf{V}\left(\mathcal{L}_{n,-\mu}\right)\times_{X\left(n\right)}S\stackrel{m_{n,\mu}}{\longrightarrow}S\]
 of the additive group $\mathbb{G}_{a,X}$ by means of $X\left(n\right)$-automorphisms
of $S$. 

\begin{prop}
Every nontrivial $\mathbb{G}_{a,X}$-action on $S$ appear in this
way.
\end{prop}
\begin{proof}
A connected component of $\pi^{-1}\left(x_{0}\right)$ is invariant
under a $\mathbb{G}_{a,X}$-action. Therefore, a $\mathbb{G}_{a,X}$-action
on $S$ lifts to a $\mathbb{G}_{a,X\left(n\right)}$-action $m:\mathbb{G}_{a,X\left(n\right)}\times_{X\left(n\right)}S\rightarrow S$.
In turn, this action restricts on $S\mid_{X_{i}\left(n\right)}\simeq\mathbf{V}\left(\mathcal{L}_{n,-\mu}\right)\mid_{X_{i}\left(n\right)}$
to a $\mathbb{G}_{a,X_{i}\left(n\right)}$-action $m_{i}$. Thus there
exists a nonzero section $s_{i}\in\Gamma\left(X_{i}\left(n\right),\mathcal{L}_{n,\mu}\right)\simeq\textrm{Hom}_{X_{i}\left(n\right)}\left(\mathcal{L}_{n,-\mu},\mathcal{O}_{X_{i}\left(n\right)}\right)$
such that the corresponding group co-action is induced by the $\mathcal{O}_{X_{i}\left(n\right)}$-modules
homomorphism $\left(\textrm{Id}\otimes1+s_{i}\otimes T\right):\mathcal{L}_{n,-\mu}\mid_{X_{i}\left(n\right)}\rightarrow\mathbf{S}\left(\mathcal{L}_{n,-\mu}\right)\mid_{X_{i}\left(n\right)}\otimes_{\mathcal{O}_{X_{i}\left(n\right)}}\mathcal{O}_{X_{i}\left(n\right)}\left[T\right].$
Clearly, these actions coincide on the overlaps $X_{ij}\left(n\right)$
if and only if the $s_{i}$'s glue to a section $s\in\Gamma\left(X,\mathcal{L}_{n,\mu}\right)$
such that $m=m_{s}$. 
\end{proof}
\begin{cor}
\label{cor:TrivialCan_equiv_freeGa} A DFS $\pi:S\rightarrow X$ admits
a free $\mathbb{G}_{a,X}$-action if and only if its canonical sheaf
is trivial.
\end{cor}
\begin{proof}
By construction, the canonical sheaf of $S=\mathbf{W}\left(\mathcal{U}_{\left(n\right)},\mathcal{L}_{n,-\mu},g\right)$
is isomorphic to $\rho_{n}^{*}\mathcal{L}_{n,-\mu}$. It is trivial
if and only if $\mathcal{L}_{n,-\mu}$ is. On the other hand, the
discussion above implies that $S$ admits a free $\mathbb{G}_{a,X}$-action
if and only if $\mathcal{L}_{n,\mu}\simeq\mathcal{L}_{n,-\mu}^{\vee}$
is trivial. 
\end{proof}

\section{Danielewski-Fieseler surfaces and labelled rooted trees}

Here we give a combinatorial description of DFS's by means of labelled
trees. To state the main result of this section, we need the following
definition. 

\begin{defn}
\emph{A relative Danielewski-Fieseler surface with base $\left(A,x\right)$}
is a morphism $\mathfrak{S}=\left(\psi:S\rightarrow\ba_{X}^{1}\right)$
of DFS's with base $\left(A,x\right)$. A \emph{morphism of relative}
\emph{DFS}'s is a morphism of DFS's $\beta:S'\rightarrow S$ such
that $\psi'=\psi\circ\beta$. 
\end{defn}
\noindent This section is devoted to the proof of the following result.

\begin{thm}
\label{thm:Equivalence_of_Cats} The category $\overrightarrow{\mathfrak{D}}_{\left(A,x\right)}$
of relative DFS's with base $\left(A,x\right)$ is equivalent to the
category $\mathcal{T}_{\left(A,x\right)}$ of $\left(A,x\right)$-labelled
rooted trees. 
\end{thm}
\noindent In the following subsections, we construct an equivalence
of categories in the form of a covariant functor $\mathfrak{S}:\mathcal{T}_{\left(A,x\right)}\rightarrow\overrightarrow{\mathfrak{D}}_{/\left(A,x\right)}$.

\subsection{DFS's defined by labelled rooted trees}

Given an $\left(A,x\right)$-labelled tree $\gamma=\left(\Gamma,\sigma\right)$
with the leaves $e_{1,m_{1}},\ldots,e_{n,m_{n}}$, we let $\mu=\mu\left(\gamma\right)=\left(\mu_{i}=h-m_{i}\right)_{i=1,\ldots,n}\in\mathbb{Z}_{\geq0}^{n}$,
where $h=h\left(\Gamma\right)$ denotes the height of $\Gamma$. Since
$\mu_{i}\leq h$ for every $i=1,\ldots,n$, the multiplication by
the regular section $x^{h}s_{-\mu}$ of $\mathcal{L}_{n,-\mu}$ defines
an homomorphism $\theta_{\mu,h}:\mathcal{L}_{n,\mu}=\mathcal{L}_{n,-\mu}^{\vee}\rightarrow\mathcal{O}_{X\left(n\right)}$
restricting to an isomorphism over $\delta_{n}^{-1}\left(X_{*}\right)$.
We let $g\left(\gamma\right)\in C^{1}\left(\mathcal{U}_{\left(n\right)},\mathcal{L}_{n,\mu}\right)$
be the unique cocycle such that $\theta_{\mu,h}\left(g\left(\gamma\right)\right)=\partial\left(\sigma\right)\in C^{1}\left(\mathcal{U}_{\left(n\right)},\mathcal{O}_{X\left(n\right)}\right)$.
If $n=1$, then the scheme $\pi_{\gamma}:S\left(\gamma\right)=\mathbf{W}\left(\mathcal{U}_{\left(n\right)},\mathcal{L}_{n,-\mu},g\left(\gamma\right)\right)\rightarrow X$
corresponding to this data is isomorphic to $\ba_{X}^{1}$. Otherwise,
if $n\geq2$, then, by (\ref{def:W_Rooted_Tree}), the transition
functions $\tilde{g}_{ij}=x^{\mu_{i}}g_{ij}=x^{-m_{i}}\left(\sigma_{j}-\sigma_{i}\right)\in A_{x}$
and $\tilde{g}_{ji}$ have a pole at $x_{0}$. Thus $S\left(\gamma\right)$
is a DFS by virtue of (\ref{txt:Affiness_criterion}). The morphism
of DFS's $\mathfrak{S}\left(\gamma\right)=\left(\psi_{\gamma}:S\left(\gamma\right)\rightarrow\ba_{X}^{1}\right)$
defined by the data $\left(\delta_{n},\theta_{\mu,h},\sigma\right)$
is called the \emph{canonical morphism associated with $\gamma$.}
The following result completes the first part of the proof of theorem
(\ref{thm:Equivalence_of_Cats}).

\begin{prop}
\label{pro:EssentiallySurjective} Every relative DFS is isomorphic
to $\mathfrak{S}\left(\gamma\right)$ for a suitable tree $\gamma$. 
\end{prop}
\begin{proof}
A morphism of DFS's $\mathfrak{S}=\left(\psi:S=\mathbf{W}\left(\mathcal{U}_{\left(n\right)},\mathcal{L}_{n,-\mu},g\right)\rightarrow\ba_{X}^{1}\right)$
is given by a data $\left(\delta_{n},\theta=a\theta_{\mu,h},\sigma\right)$,
where $a\in A^{*}$ and $h\geq\max\left(\mu_{i}\right)$, such that
$\theta\left(g\right)=\partial\left(\sigma\right)$. If $n\geq2$
then $m_{i}=h-\mu_{i}\geq1$ for every $i=1,\ldots,n$. Indeed, otherwise
there exists an indice $i$ for which $\theta$ induces an isomorphism
$\theta_{i}:\mathcal{L}_{n,\mu}\mid_{X_{i}\left(n\right)}\stackrel{\sim}{\rightarrow}\mathcal{O}_{X_{i}\left(n\right)}$.
Thus, for every $j\neq i$, $g_{ij}\in\Gamma\left(X_{ij}\left(n\right),\mathcal{L}_{n,\mu}\right)$
extends to a section of $\mathcal{L}_{n,\mu}\mid_{X_{i}\left(n\right)}$
as $\theta\left(g_{ij}\right)=\sigma_{j}\mid_{X_{i}\left(n\right)}-\sigma_{i}\mid_{X_{i}\left(n\right)}\in\Gamma\left(X_{ij}\left(n\right),\mathcal{O}_{X\left(n\right)}\right)$
extends to a section of $\mathcal{O}_{X_{i}\left(n\right)}$. This
implies that $S\mid_{X_{i}\left(n\right)\cup X_{j}\left(n\right)}$
is a trivial torsor, in contradiction with (\ref{pro:DFSCriterion}).
For the same reason, $d_{ij}=\textrm{ord}_{x_{0}}\left(\sigma_{j}-\sigma_{i}\right)<\min\left(m_{i},m_{j}\right)$
for every $i\neq j$. Since $\min\left(d_{ij},d_{ik}\right)=\min\left(d_{ij},d_{jk}\right)$
for every triple of indices $i,j$ and $k$, we deduce from (\ref{Pro:TreeConstructProp})
that the data $\left(n,\left(m_{i}\right)_{i=1,\ldots,n},\left(d_{ij}\right)_{i,j=1,\ldots,n},\sigma\right)$
corresponds to an $\left(A,x\right)$-labelled tree $\gamma=\left(\Gamma,\sigma\right)$.
Finally, the data $\left(\textrm{Id}_{X\left(n\right)},a\textrm{Id}_{\mathcal{L}_{n,-\mu}},0\right)$
defines an isomorphism of DFS's $\phi:S\stackrel{\sim}{\rightarrow}S\left(\gamma\right)$
such that $\psi_{\gamma}=\psi\circ\phi^{-1}:S\left(\gamma\right)\rightarrow\ba_{X}^{1}$. 
\end{proof}
\begin{cor}
Every DFS is $X$-isomorphic to $S\left(\gamma\right)$ for a suitable
tree $\gamma$. 
\end{cor}
\begin{proof}
Indeed, by (\ref{exa:Mor_to_A1}), every DFS $S$ admits a morphism
of DFS's $\mathfrak{S}=\left(\psi:S\rightarrow\ba_{X}^{1}\right)$. 
\end{proof}
\begin{example}
\label{exa:ODS} Given three integers $r\geq0$, $m\geq1$ and $n\geq1$,
we consider an $\left(A,x\right)$-labelled tree $\gamma=\left(\Gamma_{r,m,n},\sigma\right)$
with the following underlying tree 

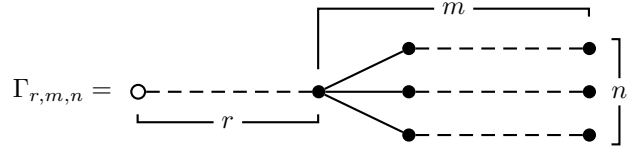
\begin{pspicture}(-2,1.8)(10,-1.3)

\def\dedge{\ncline[linestyle=dashed]}

\rput(0,0){$\Gamma_{r,m,n}=$}

\rput(4,0){

\pstree[treemode=R,radius=2.5pt,treesep=0.4cm,levelsep=1.2cm]{\Tc{3pt}}{

\skiplevels{1} 

   \pstree{\TC*[edge=\dedge]}{ 

     \pstree{\TC*}{\skiplevels{1} \TC*[edge=\dedge] \endskiplevels}

     \pstree{\TC*}{\skiplevels{1} \TC*[edge=\dedge] \endskiplevels}

    \pstree{\TC*}{\skiplevels{1}  \TC*[edge=\dedge] \endskiplevels}

   } 

\endskiplevels

}

}

\pnode(1,-0.3){A}\pnode(3.4,-0.3){B}

\ncbar[arm=3pt,angleA=-90]{A}{B}\ncput*[npos=1.5]{$r$}

\pnode(3.4,0.3){C}\pnode(7,1){D}

\ncbar[arm=3pt,angleA=90]{C}{D}\ncput*[npos=1.5]{$m$}

\pnode(7.3,0.7){E}\pnode(7.3,-0.7){F}

\ncbar[arm=3pt]{E}{F}\ncput*[npos=1.5]{$n$}

\end{pspicture}

\noindent Since there exists $b\in A$ and a cochain $\tilde{\sigma}\in A^{n}$
such that $\sigma_{i}=b+x^{r}\tilde{\sigma}_{i}$ for every $i=1,\ldots,n$,
the corresponding DFS $\pi_{\gamma}:S\left(\gamma\right)\rightarrow X$
is obtained by gluing $n$ copies $\textrm{Spec}\left(A\left[T_{i}\right]\right)$
of $\ba_{X}^{1}$ by means of the $A_{x}$-algebras isomorphisms \[
A_{x}\left[T_{i}\right]\stackrel{\sim}{\longrightarrow}A_{x}\left[T_{j}\right],\quad T_{i}\mapsto x^{-m}\left(\tilde{\sigma}_{j}-\tilde{\sigma}_{i}\right)+T_{i},\quad i\neq j.\]
Since $\tilde{\sigma}_{j}-\tilde{\sigma}_{i}\in A\setminus xA$ for
every $i\neq j$, the local sections $\tilde{\sigma}_{i}+x^{m}T_{i}\in A\left[T_{i}\right]$,
$i=1,\ldots,n$, glue to a global one $s_{1}\in B=\Gamma\left(X,\pi_{\gamma*}\mathcal{O}_{S\left(\gamma\right)}\right)$
which distinguishes the irreducible components of $\pi_{\gamma}^{-1}\left(x_{0}\right)$.
Letting $P=\prod_{i=1}^{n}\left(y-\tilde{\sigma}_{i}\right)\in A\left[y\right]$,
the rational section $x^{-m}P\left(s_{1}\right)\in B\otimes_{A}A_{x}$
extends to a regular section $s_{2}\in B$, inducing a coordinate
function on every irreducible component of $\pi_{\gamma}^{-1}\left(x_{0}\right)$.
By construction, the $A$-algebras homomorphism $A\left[y,z\right]\rightarrow B$,
$y\mapsto s_{1}$, $z\mapsto s_{2}$ defines a closed embedding $\phi:S\left(\gamma\right)\hookrightarrow\ba_{X}^{2}=\textrm{Spec}\left(A\left[y,z\right]\right)$
which induces an $X$-isomorphism between $S\left(\gamma\right)$
and the surface \[
\pi:S_{P,m}=\textrm{Spec}\left(A\left[y,z\right]/\left(x^{m}z-P\left(y\right)\right)\right)\longrightarrow X.\]
 In this coordinates, the canonical morphism $\psi_{\gamma}:S\left(\gamma\right)\rightarrow\ba_{X}^{1}$
is given as the restriction on $S_{P,m}$ of the $X$-morphism $\ba_{X}^{2}\rightarrow\ba_{X}^{1}$,
$\left(y,z\right)\mapsto x^{r}y+b$. 
\end{example}

\subsection{Morphism of DFS's defined by a morphism labelled rooted trees}

Given two $\left(A,x\right)$-labelled trees $\gamma'=\left(\Gamma',\sigma'\right)$
and $\gamma=\left(\Gamma,\sigma\right)$, we equip the corresponding
DFS's $S\left(\gamma'\right)=\mathbf{W}\left(\mathcal{U}_{\left(n'\right)},\mathcal{L}_{n',-\mu'},g'\right)$
and $S\left(\gamma\right)=\mathbf{W}\left(\mathcal{U}_{\left(n\right)},\mathcal{L}_{n,-\mu},g\right)$
with their canonical morphisms $\psi_{\gamma'}$ and $\psi_{\gamma}$
respectively. By (\ref{txt:Tree_Mor_Desc}), the image of a leaf $e'_{i,m'_{i}}$
of $\Gamma'$ by a morphism of $\left(A,x\right)$-labelled tree $\tau:\gamma'\rightarrow\gamma$
is a leaf $e_{j\left(i\right),m_{\left(j\right)}}$ of $\Gamma$ such
that $m'_{i}\geq m_{j\left(i\right)}$, and $\sigma'_{i}-\sigma_{j\left(i\right)}\in x^{m_{j\left(i\right)}}A$
for every $i=1,\ldots,n'$. Letting $\alpha:X\left(n'\right)\rightarrow X\left(n\right)$
be the unique $X$-morphism such that $\alpha\left(x_{i}\left(n'\right)\right)=x_{j\left(i\right)}\left(n\right)$
for every $i=1,\ldots,n$, we conclude that the invertible sheaf $\alpha^{*}\mathcal{L}_{n,-\mu}$
is isomorphic to $\mathcal{L}_{n',-\nu}$, where $\nu=\left(\mu_{j\left(i\right)}\right)_{i=1,\ldots,n'}\in\mathbb{Z}_{\geq0}^{n'}$.
The multiplication by the regular section $x^{h\left(\Gamma'\right)-h\left(\Gamma\right)}s_{\nu-\mu'}\in\Gamma\left(X\left(n'\right),\mathcal{L}_{n',\nu-\mu'}\right)$
defines an homomorphism of $\mathcal{O}_{X\left(n'\right)}$-modules
$\theta:\mathcal{L}_{n',\mu'}\rightarrow\alpha^{*}\mathcal{L}_{n,\mu}\simeq\mathcal{L}_{n',\nu}$
such that $\theta_{\mu',h'}=\alpha^{*}\theta_{\mu,h}\circ\theta$.
By construction, there exists a unique cochain $\sigma''\in C^{0}\left(\mathcal{U}_{\left(n'\right)},\mathcal{L}_{n',\nu}\right)$
such that $\alpha^{*}\theta_{\mu,h}\left(\sigma''\right)=\left\{ \sigma'_{i}-\sigma_{j\left(i\right)}\right\} _{i=1,\ldots n'}\in C^{0}\left(\mathcal{U}_{\left(n'\right)},\mathcal{O}_{X\left(n'\right)}\right)$.
Since $\alpha^{*}\theta_{\mu,h}$ restricts to an isomorphism over
$\delta_{n'}^{-1}\left(X_{*}\right)$, we conclude that $\theta\left(g_{\sigma'}\left(\gamma'\right)\right)=\overline{g_{\sigma}\left(\gamma\right)}+\partial\left(\sigma''\right)$,
where $\overline{g_{\sigma}\left(\gamma\right)}\in C^{1}\left(\mathcal{U}_{\left(n'\right)},\mathcal{L}_{n',\nu}\right)$
denotes the image of $\alpha^{*}g_{\sigma}\left(\gamma\right)\in C^{1}\left(\alpha^{-1}\mathcal{U}_{\left(n\right)},\alpha^{*}\mathcal{L}_{n,\mu}\right)$
by the restriction maps between \v{C}ech complexes. By (\ref{Pro:MorphismDescription}),
the data $\left(\alpha,\theta,\sigma''\right)$ defines a morphism
of DFS's $\beta_{\tau}:S\left(\gamma'\right)\rightarrow S\left(\gamma\right)$
such that $\psi_{\gamma'}=\psi_{\gamma}\circ\beta_{\tau}$, \emph{i.e.}
a morphism $\beta_{\tau}:\mathfrak{S}\left(\gamma'\right)\rightarrow\mathfrak{S}\left(\gamma\right)$
in $\overrightarrow{\mathfrak{D}}_{\left(A,x\right)}$. We say that
$\beta_{\tau}$ is the morphism of relative DFS's defined by $\tau:\gamma'\rightarrow\gamma$. 

\begin{example}
We consider the blow-down $\tau_{e}:\gamma'\rightarrow\gamma$ of
the leaves at $e$ between the following $\left(k\left[x\right],x\right)$-labelled
trees $\gamma'=\left(\Gamma',\sigma'\right)$ and $\gamma=\left(\Gamma,\sigma\right)$.

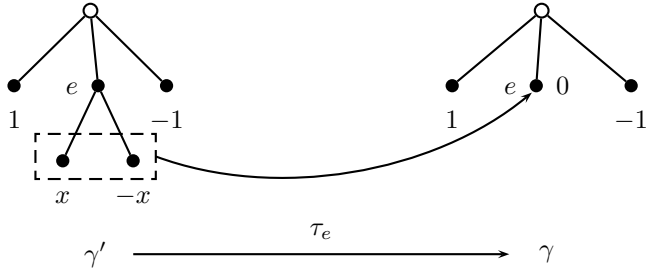
\begin{pspicture}(-2,1.8)(10,-2.5)

\def\dedge{\ncline[linestyle=dashed]}

\rput(1,0){

\pstree[treemode=D,radius=2.5pt,treesep=0.6cm,levelsep=1cm]{\Tc{3pt}}{

\TC*~[tnpos=b]{$1$} \pstree{\TC*~[tnpos=l]{$e$}}{\TC*~[tnpos=b]{$x$} \TC*~[tnpos=b]{$-x$}}\TC*~[tnpos=b]{$-1$}

}

}

\rput(7,0.5){

\pstree[treemode=D,radius=2.5pt,treesep=0.6cm,levelsep=1cm]{\Tc{3pt}}{

\TC*~[tnpos=b]{$1$} \TC*~[tnpos=l]{$e$}~[tnpos=r]{$0$} \TC*~[tnpos=b]{$-1$} 

}

}

\rput(1,-2){$\gamma'$}

\rput(7,-2){$\gamma$}

\pnode(1.5,-2){A}\pnode(6.5,-2){B}

\ncline{->}{A}{B}\naput{$\tau_e$}

\pnode(0.2,-0.7){C}\pnode(1.8,-0.7){D} \ncbox[boxsize=0.3,linestyle=dashed]{C}{D}

\pnode(6.8,0.15){F}\ncarc[arcangle=-30,angleA=0,angleB=-90, arm=10]{->}{D}{F}

\end{pspicture}

\noindent We let $s'_{1}\in B'=\Gamma\left(X,\left(\pi_{\gamma'}\right)_{*}\mathcal{O}_{S\left(\gamma'\right)}\right)$
and $s_{1}\in B=\Gamma\left(X,\left(\pi_{\gamma}\right)_{*}\mathcal{O}_{S\left(\gamma\right)}\right)$
be the sections corresponding to the canonical morphisms $\psi_{\gamma'}:S\left(\gamma'\right)\rightarrow\ba_{X}^{1}=\textrm{Spec}\left(k\left[x\right]\left[y\right]\right)$
and $\psi_{\gamma}:S\left(\gamma\right)\rightarrow\ba_{X}^{1}$ respectively.
By (\ref{exa:ODS}), $S\left(\gamma\right)$ is isomorphic to the
surface $S=\left\{ xz-y\left(y^{2}-1\right)=0\right\} \subset\ba_{X}^{2}=\textrm{Spec}\left(k\left[x\right]\left[y,z\right]\right)$
via the embedding induced by the $k\left[x\right]$-algebra homomorphism
$k\left[x\right]\left[y,z\right]\mapsto B$, $\left(y,z\right)\mapsto\left(s_{1},s_{2}=x^{-1}s_{1}\left(s_{1}^{2}-1\right)\right)$.
A similar argument shows that $S\left(\gamma'\right)$ is isomorphic
to the Bandman and Makar-Limanov surface $S'\subset\ba_{X}^{3}=\textrm{Spec}\left(k\left[x\right]\left[y,z,u\right]\right)$
with equations \[
xz-y\left(y^{2}-1\right)=0,\quad yu-z\left(z^{2}-1\right)=0,\quad xu-\left(y^{2}-1\right)\left(z^{2}-1\right)=0,\]
 via the embedding induced by the $k\left[x\right]$-algebras homomorphism
$k\left[x\right]\left[y,z,u\right]\rightarrow B$ \[
\left(y,z,u\right)\mapsto\left(s'_{1},s'_{2}=x^{-1}s'_{1}\left(\left(s'_{1}\right)^{2}-1\right),s'_{3}=x^{-1}\left(\left(s'_{1}\right)^{2}-1\right)\left(\left(s'_{2}\right)^{2}-1\right)\right)\]
In these coordinates, the morphism $\beta_{\tau_{e}}:S\left(\gamma'\right)\rightarrow S\left(\gamma\right)$
defined by the blow-down $\tau_{e}$ coincide with the restriction
of the projection $\ba_{X}^{3}\rightarrow\ba_{X}^{2}$, $\left(x,y,z,u\right)\mapsto\left(x,y,z\right)$. 
\end{example}
\begin{enavant} The correspondence $\gamma\mapsto\mathfrak{S}\left(\gamma\right)$,
$\left(\tau:\gamma'\rightarrow\gamma\right)\mapsto\mathfrak{S}\left(\tau\right)=\left(\beta_{\tau}:\mathfrak{S}\left(\gamma'\right)\rightarrow\mathfrak{S}\left(\gamma\right)\right)$
defines a covariant functor $\mathfrak{S}:\mathcal{T}_{\left(A,x\right)}\rightarrow\overrightarrow{\mathfrak{D}}_{\left(A,x\right)}$.
It follows from (\ref{pro:EssentiallySurjective}) that $\mathfrak{S}$
is essentially surjective. The following result shows that $\mathfrak{S}$
is fully faithful, whence completes the proof of theorem (\ref{thm:Equivalence_of_Cats}). 

\end{enavant}

\begin{prop}
\label{pro:FullyFaithful} Every morphism $\beta:\mathfrak{S}\left(\gamma'\right)\rightarrow\mathfrak{S}\left(\gamma\right)$
coincides with a unique morphism $\beta_{\tau}$ defined by a morphism
of $\left(A,x\right)$-labelled trees $\tau:\gamma'\rightarrow\gamma$. 
\end{prop}
\begin{proof}
By (\ref{Pro:MorphismDescription}), $\beta$ is determined by a data
$\left(\alpha,\theta,\sigma''\right)$ such that $\alpha\left(x_{i}\left(n'\right)\right)=x_{j\left(i\right)}\left(n\right)$,
$i=1,\ldots,n'$, $\theta_{\mu',h}=\alpha^{*}\theta_{\mu,h}\circ\theta$,
$\sigma'=\alpha^{*}\sigma+\alpha^{*}\theta_{\mu,h}\left(\sigma''\right)$,
and such that \[
\theta\left(g_{\sigma'}\left(\gamma'\right)\right)=\overline{\alpha^{*}g_{\sigma}\left(\gamma\right)}+\partial\left(\sigma''\right)\in C^{1}\left(\mathcal{U}_{\left(n'\right)},\alpha^{*}\mathcal{L}_{n,\mu}\right).\]
 Since $\alpha^{*}\mathcal{L}_{n,\mu}$ is isomorphic to $\mathcal{L}_{n',\nu}$,
where $\nu=\left(\mu_{j\left(i\right)}\right)_{i}\in\mathbb{Z}_{\geq0}^{n'}$,
we conclude that $\theta$ is the multiplication by $x^{h'-h}s_{\nu-\mu'}\in\Gamma\left(X\left(n'\right),\mathcal{L}_{n',\nu-\mu'}\right)$.
Thus $m_{i}'\geq m_{j\left(i\right)}$ and $\tilde{\sigma}_{i}''=\alpha^{*}\theta_{\mu,h}\left(\sigma_{i}''\right)\in x^{h-\mu_{j\left(i\right)}}A=x^{m_{j\left(i\right)}}A$
for every $i=1,\ldots,n'$. Therefore, the formulas \[
\Gamma'_{i}\ni e'_{i,k}\mapsto\tau_{i}\left(e'_{i,k}\right)=e_{j\left(i\right),\min\left(k,m_{j\left(i\right)}\right)}\in\Gamma_{j\left(i\right)},\quad i=1,\ldots,n'\; k=0,\ldots,m'_{i},\]
(see (\ref{eq:MaxChains}) for the notation) define a collection of
morphisms of trees $\tau_{i}:\left(\Gamma'_{i},\sigma'_{i}\right)\rightarrow\left(\Gamma_{j\left(i\right)},\sigma_{j\left(i\right)}\right)$
in $\mathcal{T}_{\left(A,x\right)}$. Moreover, we deduce that \begin{equation}
\textrm{ord}_{x_{0}}\left(\sigma'_{i'}-\sigma'_{i}\right)=\left\{ \begin{array}{ccc}
\textrm{ord}_{x_{0}}\left(\tilde{\sigma}''_{i'}-\tilde{\sigma''}_{i}\right)\geq m_{j\left(i\right)} &  & \textrm{if }j\left(i\right)=j\left(i'\right)\\
\textrm{ord}_{x_{0}}\left(\sigma_{j\left(i'\right)}-\sigma_{j\left(i\right)}\right) &  & \textrm{otherwise}\end{array}\right..\label{SameOrder}\end{equation}
 Since $\sigma'$ and $\sigma$ are compatible with $\Gamma'$ and
$\Gamma$ respectively, (\ref{SameOrder}) guarantees that the conditions
of (\ref{Pro:MorphismConstructionProp}) are satisfied. Thus there
exists a unique morphism $\tau:\gamma'\rightarrow\gamma$ in $\mathcal{T}_{\left(A,x\right)}$
such that $\tau\mid_{\gamma_{i}'}=\tau_{i}$ for every $i=1,\ldots,n'$.
By construction, $\beta=\beta_{\tau}$. 
\end{proof}
\begin{cor}
For every morphism of DFS's $\beta:S'\rightarrow S$ there exist a
morphism of $\left(A,x\right)$-labelled trees $\tau:\gamma'\rightarrow\gamma$
and a commutative diagram \[\xymatrix{ S' \ar[r]^{\beta} \ar[d]^{\wr}_{\phi'} & S \ar[d]^{\phi}_{\wr} \\ S\left(\gamma'\right) \ar[r]^{\beta_{\tau}} & S\left(\gamma\right) }\] 
\end{cor}
\begin{proof}
Once a morphism of DFS's $\psi:S\rightarrow\ba_{X}^{1}$ is chosen
(see (\ref{exa:Mor_to_A1})), $\left(\psi:S\rightarrow\ba_{X}^{1}\right)$
and $\left(\psi'=\psi\circ\beta:S'\rightarrow\ba_{X}^{1}\right)$
are objects of $\overrightarrow{\mathfrak{D}}_{\left(A,x\right)}$,
whereas $\beta:S'\rightarrow S$ corresponds to morphism in $\overrightarrow{\mathfrak{D}}_{\left(A,x\right)}$.
So the result follows from (\ref{pro:EssentiallySurjective}) and
(\ref{pro:FullyFaithful}). 
\end{proof}

\subsection{Reading isomorphism classes of Danielewski-Fieseler surfaces from
trees}

To decide when two DFS's are $X$-isomorphic, we have the following
criterion.

\begin{thm}
\label{EquivalentTreesTheorem} Two $\left(A,x\right)$-labelled trees
define $X$-isomorphic DFS's if and only if they are equivalent (see
(\ref{def:Equivalence_of_trees}) for the definition). 
\end{thm}
\begin{proof}
We first observe that if $\gamma'=\left(\Gamma',\sigma'\right)$ is
an essential labelled subtree for $\gamma=\left(\Gamma,\sigma\right)$
then the DFS's $S\left(\gamma\right)$ and $S\left(\gamma'\right)$
are $X$-isomorphic. Indeed, there exists $m\in\mathbb{Z}_{\geq0}$
and $b\in A$ such that $m_{i}=m'_{i}+m$ and $\sigma_{i}=b+x^{m}\sigma'_{i}$
for every $i=1,\ldots,n$. Thus $\mu\left(\gamma'\right)=\mu\left(\gamma\right)=\mu\in\mathbb{Z}_{\geq0}^{n}$,
$g\left(\gamma\right)=g\left(\gamma'\right)\in C^{1}\left(\mathcal{U}_{\left(n\right)},\mathcal{L}_{n,\mu}\right)$,
and so, $S\left(\gamma\right)\simeq S\left(\gamma'\right)$. Therefore,
it suffices to prove the assertion for DFS's defined by essential
labelled trees $\gamma=\left(\Gamma,\sigma\right)$ and $\gamma'=\left(\Gamma',\sigma'\right)$.
If $\gamma$ and $\gamma'$ are equivalent then $h=h\left(\Gamma\right)=h\left(\Gamma'\right)$
and there exists a permutation $j$ of $\left\{ 1,\ldots,n\right\} $
such that $\mu'_{i}=\mu_{j\left(i\right)}$ for every $i=1,\ldots,n$.
Moreover, there exist a pair $\left(a,b\right)\in A^{*}\times A$
and a cochain $\tilde{\sigma}=\left\{ \tilde{\sigma}_{i}\right\} _{i=1,\ldots,n}\in A^{n}$
such that $\sigma'_{i}=a\sigma_{j\left(i\right)}+b+x^{m_{j\left(i\right)}}\tilde{\sigma}_{j\left(i\right)}$
for every $i=1,\ldots,n$. Letting $\sigma''=\left\{ \sigma_{i}''=\sigma'_{i}-b\right\} _{i=1,\ldots,n}$,
the same argument as above shows that $S\left(\gamma'\right)$ is
$X$-isomorphic to the DFS defined by the labelled tree $\gamma''=\left(\Gamma',\sigma''\right)$.
Thus, by replacing $S\left(\gamma'\right)$ by $S\left(\gamma''\right)$
if necessary, we can suppose from now on that $b=0$. We let $\alpha\in\textrm{Aut}_{X}\left(X\left(n\right)\right)$
be the unique $X$-automorphism such that $\alpha\left(x_{i}\left(n\right)\right)=x_{j\left(i\right)}\left(n\right)$,
and we let $c\in C^{0}\left(\mathcal{U}_{\left(n\right)},\alpha^{*}\mathcal{L}_{n,\mu}\right)$
be the unique cochain such that $\theta_{\mu',h}\left(c\right)=\left\{ a^{-1}x^{m_{j\left(i\right)}}\tilde{\sigma}_{j\left(i\right)}\right\} _{i=1,\ldots,n}$.
Since $\mathcal{L}_{n,\mu'}\simeq\alpha^{*}\mathcal{L}_{n,\mu}$ and
$a^{-1}g\left(\gamma'\right)=\alpha^{*}g\left(\gamma\right)+\partial\left(c\right)$,
we conclude that the data $\left(\alpha,a\textrm{Id}_{\mathcal{L}_{n,\mu'}},c\right)$
defines an isomorphism of DFS's $\beta:S\left(\gamma'\right)\stackrel{\sim}{\rightarrow}S\left(\gamma\right)$.
Conversely, an $X$-isomorphism $\beta:S'=S\left(\gamma'\right)\stackrel{\sim}{\rightarrow}S=S\left(\gamma\right)$
exists if and only if the trees $\gamma$ and $\gamma'$ have the
same number $n$ of leaves, and there exists a data \[
\left(\alpha,\theta,c\right)\in\textrm{Aut}_{X}\left(X\left(n\right)\right)\times\textrm{Isom}_{X\left(n\right)}\left(\mathcal{L}_{n,\mu'},\alpha^{*}\mathcal{L}_{n,\mu}\right)\times C^{0}\left(\mathcal{U}_{\left(n\right)},\alpha^{*}\mathcal{L}_{n,\mu}\right)\]
 such that $\theta\left(g\left(\gamma'\right)\right)=\alpha^{*}g\left(\gamma\right)+\partial\left(c\right)$.
By replacing $S$ by the $X$-isomorphic DFS $\alpha^{*}S$, we can
suppose that $\alpha=\textrm{Id}_{X\left(n\right)}$. We can also
suppose that $h'=h\left(\Gamma'\right)\geq h=h\left(\Gamma\right)$.
Then $\mathcal{L}_{n,\mu'}$ and $\mathcal{L}_{n,\mu}$ are isomorphic
if and only if $m'_{i}=h'-\mu'_{i}=h-\mu_{i}+m=m_{i}+m$ for some
$m\in\mathbb{Z}_{\geq0}$. The identity $\theta\left(g\left(\gamma'\right)\right)=g\left(\gamma\right)+\partial\left(c\right)$
is satisfied if and only there exist a pair $\left(a,b\right)\in A^{*}\times A$
such that $\tilde{\sigma}_{i}:=\sigma_{i}+\theta_{\mu,h}\left(c_{i}\right)=ax^{m}\sigma'_{i}+b\in A$
for every $i=1,\ldots,n$. By definition, $\theta_{\mu,h}\left(c_{i}\right)\in x^{m_{i}}A\subset xA$
for every $i=1,\ldots,n$. Moreover, since $\Gamma$ is essential,
there exist $i\neq j$ such that $\sigma_{j}-\sigma_{i}\in A\setminus xA$.
Thus $m=0$, and so $\textrm{ord}_{x_{0}}\left(\sigma_{j}-\sigma_{i}\right)=\textrm{ord}_{x_{0}}\left(\sigma'_{j}-\sigma'_{i}\right)$
for every $i\neq j$. Letting $\tilde{\sigma}=\left\{ \tilde{\sigma}_{i}\right\} _{i=1,\ldots,n}$,
we conclude that the trees $\gamma'$ and $\tilde{\gamma}=\left(\Gamma,\tilde{\sigma}\right)$
are equivalent. This completes the proof as $\gamma$ and $\tilde{\gamma}$
are isomorphic, whence equivalent.
\end{proof}
\begin{example}
The following $\left(k\left[x\right],x\right)$-labelled trees $\gamma_{t}=\left(\Gamma,\sigma_{t}\right)$,
where $t\in k$, 

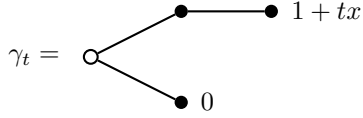
\begin{pspicture}(-2.2,-0.8)(8,1.2)

\rput(1.5,-0.02){$\gamma_t=$}

\rput(4,0){

\pstree[treemode=R,radius=2.5pt,treesep=1cm,levelsep=1.2cm]{\Tc{3pt}}{

\pstree{\TC*}{\TC*~[tnpos=r]{$1+tx$}} \TC*~[tnpos=r]{$0$}

}

}

\end{pspicture}

\noindent define a family of DFS's $\pi_{t}:S_{t}=S\left(\gamma_{t}\right)\rightarrow X=\textrm{Spec}\left(k\left[x\right]\right)$.
A similar argument as in (\ref{exa:ODS}) shows that $S_{t}$ is $X$-isomorphic
to the fiber $pr_{v}^{-1}\left(t\right)$ of the threefold $pr_{v}:\mathcal{X}\subset\textrm{Spec}\left(k\left[x\right]\left[y,z,u,v\right]\right)\rightarrow\textrm{Spec}\left(k\left[v\right]\right)$
with equations \[
xz-y\left(y-1\right)=0,\quad yu-z\left(z-v\right)=0,\quad xu-y\left(z-v\right)=0.\]
 These DFS's $S_{t}$ are two by two non $X$-isomorphic as the trees
$\gamma_{t}$ are two by two nonequivalent. However, there are all
isomorphic as abstract schemes. Indeed, a similar argument as in (\ref{CompletionProcedure})
below shows that there exists a projective surface $\bar{S}_{t}$,
obtained from the Hirzebruch surface $\bar{p}_{1}:\mathbb{F}_{1}=\mathbb{P}\left(\mathcal{O}_{\bp^{1}}\oplus\mathcal{O}_{\bp^{1}}\left(-1\right)\right)\rightarrow\bp_{k}^{1}$
by first blowing two distinct points $y_{1},y_{2}\in F_{0}=\bar{p}_{1}^{-1}\left(0\right)$
with exceptional divisors $E$ and $F$, and then blowing-up a point
$y_{3}\in E\setminus F_{0}'$, and an open embedding of $S_{t}$ in
$\bar{S}_{t}$ as the complement of the union of the strict transforms
of $E$, $F_{0}$, $F_{\infty}=\bar{p}_{1}^{-1}\left(\infty\right)$
and of a section $C$ of $\bar{p}_{1}$ 'at infinity'. By successively
blowing-down the strict transforms of $C$, $F_{0}$ and $E$, we
realize $S_{t}$ as the complement in $\mathbb{F}_{1}$ of a section
$D_{t}$ of $\bar{p}_{1}$ with self-intersection $\left(D_{t}^{2}\right)=3$.
By a result of Gizatullin-Danilov \cite{GizDan77}, the isomorphism
class of $\mathbb{F}_{1}\setminus D$ as an abstract scheme does not
depend on the choice of a section $D$ such that $\left(D^{2}\right)=3$.
Therefore, the surfaces $S_{t}$, $t\in k$, are all isomorphic. Since
they are not $X$-isomorphic, we deduce that a surface $S\simeq\mathbb{F}_{1}\setminus D$
as above comes equipped with a family $q_{t}:S\rightarrow X_{t}$
of structures of DFS over different bases $X_{t}\simeq\ba_{k}^{1}$,
$t\in k$, such that the general fibers of $q_{t}$ and $q_{t'}$
do not coincide whenever $t\neq t'$. 
\end{example}

\section{Danielewski-Fieseler surfaces and affine modifications}

Here we exploit the structure of rooted trees to give another description
of morphisms of DFS's. As a consequence, we obtain a canonical procedure
to construct an open embedding of a DFS $\pi:S\rightarrow X$ into
a projective $X$-scheme $\bar{\pi}:\bar{S}\rightarrow X$.

\subsection{Morphisms of Danielewski-Fieseler surfaces as affine modifications}

In what follows, we freely use the results about affine modifications,
referring the reader to \cite{KaZa99} for complete proofs. However,
we recall the following definition. 

\begin{defn}
Suppose we are given an affine scheme $V=\textrm{Spec}\left(B\right)$,
a nonzero divisor $f\in B$ an an ideal $I\subset B$ containing $f$.
The \emph{proper transform} $D_{f}^{pr}$ of $\textrm{div}\left(f\right)$
in $\bar{V}=\textrm{Proj}\left(B\left[It\right]\right)$ is the set
of prime ideals $\mathfrak{p}\in\bar{V}$ such that $ft\in\mathfrak{p}$.
The affine scheme $V'=\bar{V}\setminus\textrm{Supp}\left(D_{f}^{pr}\right)\simeq\textrm{Spec}\left(B\left[It\right]/\left(1-ft\right)\right)$
is called the \emph{affine modification} \emph{of} $V$ \emph{with
the locus} $\left(I,f\right)$. 
\end{defn}
\begin{defn} A morphism of DFS's $\beta:S'\rightarrow S$ is called
a \emph{fibered modification} if there exist an $\left(A,x\right)$-labelled
tree $\gamma'=\left(\Gamma',\sigma'\right)$, a blow-down $\tau_{e'}:\gamma'\rightarrow\gamma$
of the leaves at a certain $e'\in\Gamma'$ (see (\ref{def:TreeBlowDef})
for the definition) and a commutative diagram 

\[ \xymatrix{ S' \ar[rr]^{\beta} \ar[d]_{\wr} && S \ar[d]^{\wr} \\ S\left(\gamma'\right) \ar[rr]^{\beta_{\tau_e'}} & & S\left(\gamma\right) }\]

\end{defn}

\begin{enavant} Consider a blow-down $\tau_{e'}:\gamma'\rightarrow\gamma$
of the leaves $e_{1}',\ldots,e_{r}'$ at $e'\in\Gamma'$. For every
$e_{i}'\in\textrm{Leaves}\left(\Gamma\right)\setminus\left\{ e_{1}',\ldots,e_{r}'\right\} $,
$e_{j\left(i\right)}=\tau_{e'}\left(e_{i}'\right)$ is a leaf of $\gamma'$,
and $\tau_{e'}$ restricts to an isomorphism of labelled chains $\left(\downarrow e_{i}'\right)_{\gamma'}\stackrel{\sim}{\rightarrow}\left(\downarrow e_{j\left(i\right)}\right)_{\gamma}$.
Therefore, the morphism $\beta_{\tau_{e'}}:S'=S\left(\gamma'\right)\rightarrow S=S\left(\gamma\right)$
induces isomorphisms $S'\mid_{X_{i}\left(n'\right)}\stackrel{\sim}{\rightarrow}S\mid_{X_{j\left(i\right)}\left(n\right)}$
between the corresponding open subsets of $S'$ and $S$ respectively.
On the other hand, the image of the labelled subtree $\gamma\left(e'\right)=\left(\downarrow e'\right)_{\gamma'}\cup\textrm{Ch}_{\Gamma'}\left(e'\right)$
of $\gamma'$ by $\tau_{e'}$ is the subtree $\left(\downarrow_{\tau_{e'}\left(e'\right)}\right)_{\gamma}$
of $\gamma$. By (\ref{exa:ODS}), there exists an $X$-isomorphism
between $\tilde{S}=S'\mid_{X_{1}\left(n'\right)\cup\cdots\cup X_{r}\left(n'\right)}$
and a DFS $S_{P}=\textrm{Spec}\left(A\left[y,z\right]/\left(xz-P\left(y\right)\right)\right)$
for a certain polynomial $P\in A\left[y\right]$ of degree $r$ whose
residue class $\bar{P}\in A\left[y\right]/xA\left[y\right]\simeq k\left[y\right]$
has $r$ simple roots. In this coordinates, the restriction $\beta:\tilde{S}\rightarrow S\mid_{X_{j\left(i\right)}\left(n\right)}\simeq\ba_{X}^{1}$
coincides with the first projection $p_{y}:S_{P}\rightarrow\ba_{X}^{1}=\textrm{Spec}\left(A\left[y\right]\right)$.
Therefore, $\tilde{S}$ is isomorphic to the affine modification of
$\ba_{X}^{1}$ with locus $\left(I,f\right)=\left(\left(x,P\left(y\right)\right),x\right)$.
Since $\left(x,P\left(y\right)\right)$ is a regular sequence in $A\left[y\right]$,
the Rees algebra $A\left[y\right]\left[It\right]$ is isomorphic to
$A\left[y\right]\left[u,v\right]/\left(xv-P\left(y\right)u\right)$
via the map $u\mapsto xt$, $v\mapsto P\left(y\right)t$. Therefore,
the proper transform $D_{x}^{pr}=\left\{ xt=u=0\right\} $ of the
line $L=\textrm{div}\left(x\right)\subset\ba_{X}^{1}$ by the blow-up
morphism $\bar{\beta}:\bar{S}_{P}=\textrm{Proj}_{A}\left(A\left[y\right]\left[It\right]\right)\rightarrow\ba_{X}^{1}$
coincides with the usual strict transform $L'$ of $L$. This leads
to the following description. 

\end{enavant}

\begin{prop}
Given a fibered modification $\beta:S'\rightarrow S$, there exists
a unique irreducible component $C$ of $\pi^{-1}\left(x_{0}\right)$
such that $\beta\left(\beta^{-1}\left(C\right)\right)$ is supported
at a $0$-dimensional closed subscheme $Y\subset C$. Letting $\bar{\beta}:\bar{S}'\rightarrow S$
be the blow-up of $Y$ equipped with its reduced structure, $S'$
is isomorphic to the complement in $\bar{S}'$ of the strict transform
of $C$, and the irreducible components of $\left(\pi'\right)^{-1}\left(x_{0}\right)$
contained in $\beta^{-1}\left(C\right)$ coincides with the intersections
of the exceptional divisors of $\bar{\beta}$ with $S'$.
\end{prop}
\noindent As a consequence of (\ref{Pro:TreeBlowUpLemma}) and (\ref{thm:Equivalence_of_Cats}),
we obtain the following result.

\begin{thm}
\label{FiberedModifThm} Every morphism of Danielewski-Fieseler surfaces
factors into a finite sequence of fibered modifications followed by
an open embedding. 
\end{thm}

\subsection{Local completions of Danielewski-Fieseler surfaces}

In this subsection, we construct open embeddings of a DFS $\pi:S\rightarrow X$
into projective $X$-schemes $\bar{\pi}:\bar{S}\rightarrow X$. Our
procedure coincides with the one of Theorem 2 in \cite{Fie94} in
case that $\pi:S\rightarrow X$ is defined over $\bc$. 

\begin{enavant}\label{CompletionProcedure} Given a DFS $\pi:S=S\left(\gamma\right)\rightarrow X$
defined by an essential $\left(A,x\right)$-labelled tree $\gamma$,
we consider a factorization of the canonical morphism $\psi_{\gamma}:S=S_{m}\rightarrow\ba_{X}^{1}=S_{0}$
into a sequence of fibered modifications $\beta_{k}:S_{k}\rightarrow S_{k-1}$,
$k=1,\ldots,m$. We embed $S_{0}$ in the $\bp^{1}$-bundle $\bar{\pi}_{0}:\bar{S}_{0}:=\bp_{X}^{1}\rightarrow X$
as the complement a section $C\simeq X$ 'at infinity'. For every
$k=1,\ldots,m$, we let $\bar{\beta}_{k}:\bar{S}_{k}\rightarrow\bar{S}_{k-1}$
the blow-up morphism corresponding to the fibered modification $\beta_{k}:S_{k}\rightarrow S_{k-1}$.
By construction, the morphism $\psi_{\gamma}$ lifts to an open embedding
of $S$ into $\bar{\pi}:\bar{S}=\bar{S}_{m}\rightarrow X$. We denote
by $E_{j}\simeq\bp_{k}^{1}$, $j\in J$, the irreducible components
of $\bar{\pi}^{-1}\left(F_{0}\right)$ different from the strict transform
$F'_{0}\simeq\bp_{k}^{1}$ of $F_{0}=\bar{p}_{0}^{-1}\left(x_{0}\right)$
and the closures $\bar{C}_{e}\simeq\bp_{k}^{1}$, $e\in L\left(\Gamma\right)$,
of the irreducible components of $\pi^{-1}\left(x_{0}\right)$. Since
for every $k=1,\ldots,m$, the center of the blow-up $\bar{\beta}_{k}:\bar{S}_{k}\rightarrow\bar{S}_{k-1}$
is contained in the closure $\bar{D}_{k-1}\simeq\bp_{k}^{1}$ in $\bar{S}_{k-1}$
of a unique irreducible component of $\bar{\beta}_{k-1}^{-1}\left(F_{0}\right)\cap S_{k-1}$,
we conclude that $B_{k}=\bar{S}_{k}\setminus S_{k}$ is the strict
transform of $B_{k-1}\cup\bar{D}_{k-1}$. Therefore, the dual graph
of the $SNC$-divisor \[
\textrm{Supp}\left(\bar{\pi}^{-1}\left(x_{0}\right)\right)=\overline{\left(B_{m}\setminus C'\right)}\cup{\displaystyle \bigcup_{e\in L\left(\Gamma\right)}\bar{C}_{e}}\]
 coincides with the underlying tree $\Gamma$ of $\gamma$. Since
the tree $\Gamma$ is essential by assumption, we conclude that $\left({F'}^2_0 \right)=-\textrm{Card}\left(\textrm{Child}_{\Gamma} \left(e_0\right)\right)\neq -1$ ,
whereas $\left(E_{j}^{2}\right)\leq-2$ for every curve $E_{j}\subset\bar{S}$,
$j\in J$. This leads to the following result. 

\end{enavant}

\begin{prop}
\label{pro:CanonicalCompletion} A DFS $\pi:S=S\left(\gamma\right)\rightarrow X$
with base $\left(X,x_{0}\right)$ defined by a labelled tree $\gamma=\left(\Gamma,\sigma\right)$
admits an open embedding $i:S\hookrightarrow\bar{S}$ into a regular
projective $X$-scheme $\bar{\pi}:\bar{S}\rightarrow X$ with the
following properties.

a) $\bar{\pi}$ restricts to a trivial $\bp^{1}$-bundle over $X_{*}=X\setminus\left\{ x_{0}\right\} $. 

b) The fiber $\bar{S}_{x_{0}}=\bar{\pi}^{-1}\left(x_{0}\right)$ is
a reduced $SNC$-divisor which does not contain $\left(-1\right)$-curves
meeting at most two other components transversally in a single point. 

c) The irreducible components of $\bar{S}_{x_{0}}$ are isomorphic
to $\bp_{k}^{1}$.

d) The dual graph of $\bar{S}_{x_{0}}$ and $\bar{S}_{x_{0}}\setminus i\left(S_{x_{0}}\right)$
are isomorphic to the underlying tree $\Gamma'$ of an essential tree
for $\Gamma$ and to the subtree subtree $\Gamma'\setminus L\left(\Gamma'\right)$
of $\Gamma'$ respectively. 
\end{prop}

\section{Danielewski surfaces}

We recall (\ref{def:DFS_def}) that a \emph{Danielewski surface} \emph{is}
a DFS with base $\left(k\left[x\right],x\right)$, where $k$ denotes
a field of caracteristic zero. As a consequence of (\ref{exa:Weighted_Rooted_Trees})
and (\ref{thm:Equivalence_of_Cats}), the category $\overrightarrow{\mathfrak{D}}_{\left(k\left[x\right],x\right)}$
of relative Danielewski surfaces is equivalent to the category $\mathcal{T}_{w}^{k}$
of fine $k$-weighted trees. Therefore, every isomorphism class in
$\overrightarrow{\mathfrak{D}}_{\left(k\left[x\right],x\right)}$
contains a canonical object $\left(\psi_{0}:S_{0}\rightarrow\ba_{X}^{1}\right)$,
consisting of the canonical morphism of the DFS $\pi:S_{0}=S\left(\gamma_{w}\right)\rightarrow X$
defined by the $\left(k\left[x\right],x\right)$-labelled tree $\gamma_{w}$
obtain from a fine $k$-weighted tree by the procedure described in
(\ref{exa:Weighted_Rooted_Trees}). In this section, we characterize
Danielewski surfaces with a trivial Makar-Limanov invariant in case
that the base field $k$ is algebraically closed. 

\begin{notation}
Throughout this section, $\bar{k}$ denotes an algebraically closed
field of caracteristic zero.
\end{notation}

\subsection{Normal affine surfaces with a trivial Makar-Limanov invariant}

The Makar-Limanov invariant \cite{KML97} of an affine algebraic variety
$V=\textrm{Spec}\left(B\right)$ over $\bar{k}$ is the sub-algebra
$\textrm{ML}\left(V\right)$ of $B$ consisting of regular functions
on $V$ which are invariant under all $\mathbb{G}_{a,k}$-action on
$V$. If $\textrm{ML}\left(V\right)=\bar{k}$, then we say that $V$
has a trivial Makar-Limanov invariant. It is known that a normal affine
surface surface $S$ has a trivial Makar-Limanov invariant if and
only if it admits two $\ba^{1}$-fibrations $\pi_{1}:S\rightarrow X_{1}\simeq\ba_{\bar{k}}^{1}$
and $\pi_{2}:S\rightarrow X_{2}\simeq\ba_{\bar{k}}^{1}$ with distinct
general fibers. This property leads a geometrical criterion for $S$
to have trivial Makar-Limanov invariant. We recall that a minimal
completion of $S$ is an open embedding $S\hookrightarrow\bar{S}$
of $S$ into a normal projective surface $\bar{S}$, nonsingular along
$B=\bar{S}\setminus S$, such that $B$ is an $SNC$-divisor, with
nonsingular curves as irreducible components, which does not contain
$\left(-1\right)$-curves meeting at most two other components transversally
in a single point. In \cite{Dub02}, the author established that a
complex normal affine surface $S$ has a trivial Makar-Limanov invariant
if and only if the dual graph of the boundary divisor of any minimal
completion of $S$ is a chain. Since the arguments given in Proposition
2.10 and Theorem 2.16 in \cite{Dub02} remain valid over an arbitrary
algebraically closed field of caracteristic zero, we obtain the following
more general criterion. 

\begin{thm}
\label{thm:ML_general_carac} A normal affine surface $S\not\simeq\textrm{Spec}\left(\bar{k}\left[x,x^{-1},y\right]\right)$
has a trivial Makar-Limanov invariant if and only if the dual graph
of the boundary divisor of every minimal completion of $S$ is a chain.
Moreover, if there exists a completion of $S$ with this property,
then this also holds for every other minimal completion. 
\end{thm}

\subsection{Danielewski surfaces with a trivial Makar-Limanov invariant}

Here we apply the general criterion (\ref{thm:ML_general_carac})
to Danielewski surfaces.

\begin{defn}
\label{CombDef} A \emph{comb} is a rooted tree $\Gamma$ such that
$\Gamma\setminus L\left(\Gamma\right)$ is a chain. Equivalently,
$\Gamma$ is a comb if every $e\in\Gamma\setminus L\left(\Gamma\right)$
has at most one child which is not a leaf of $\Gamma$. For instance,
the tree $\Gamma$ of Figure $1$ is a comb rooted in $e_{0}$. 
\end{defn}
\noindent For Danielewski surfaces

\begin{thm}
\label{MLGDSThm} A Danielewski surface $S/\bar{k}$ has a trivial
Makar-Limanov invariant if and only if it is isomorphic to a Danielewski
surface $S\left(\gamma\right)$ defined by an $\left(\bar{k}\left[x\right],x\right)$-labelled
comb $\gamma$. 
\end{thm}
\begin{proof}
A tree $\Gamma$ is a comb if and only if its essential subtree $\textrm{Es}\left(\Gamma\right)$
is. Thus, by (\ref{EquivalentTreesTheorem}), we can suppose that
$S=S\left(\gamma\right)$ for a certain essential $\left(\bar{k}\left[x\right],x\right)$-labelled
tree $\gamma=\left(\Gamma,\sigma\right)$. If $\gamma$ is the trivial
tree, then $S\simeq\ba_{X}^{1}\simeq\ba_{k}^{2}$ has trivial Makar-Limanov
invariant. Otherwise we let $p_{0}:\bar{S}_{0}=\bp_{\bar{X}}^{1}\rightarrow\bar{X}$,
where $\bar{X}\simeq\bp_{X}^{1}$ denotes a nonsingular projective
model of $X\simeq\ba_{\bar{k}}^{1}$, and we embed $S_{0}=\ba_{X}^{1}$
into $\bar{S}_{0}$ as the complement of the ample divisor $F_{\infty}\cup C$,
where $F_{\infty}\simeq\bp_{\bar{k}}^{1}$ denotes the fiber of $p_{0}$
over $\infty=\bar{X}\setminus X$, and where $C\simeq\bp_{\bar{k}}^{1}$
is a section of $p_{0}$ 'at infinity'. By applying the same procedure
as in (\ref{CompletionProcedure}) to $\bar{S}_{0}$, we obtain a
projective surface $\bar{S}$ and an open embedding $S\hookrightarrow\bar{S}$
such that the dual graph $G\left(B\right)$ of the boundary divisor
$B=\bar{S}\setminus S$ is isomorphic to the tree obtained from $\Gamma$
by deleting its leaves and replacing its root by a chain with three
elements corresponding to the strict transforms $F'_{0},C'$ and $F'_{\infty}$
of the curves $F_{0}=p_{0}^{-1}\left(x_{0}\right)$, $C$ and $F_{\infty}$.
Since $\left({C'}^{2}\right)=\left({F_{\infty}'}^{2}\right)=0$ and
$\left({F_{0}'}^{2}\right)\leq-2$ as $\gamma$ is essential, we conclude
that $\bar{S}$ is a minimal completion of $S$. So the statement
follows from (\ref{thm:ML_general_carac}) as $G\left(B\right)$ is
a chain if and only if $\Gamma$ is a comb. 
\end{proof}
\begin{example}
Let $P_{1},\ldots,P_{n}\in\bar{k}\left[T\right]$ be a collection
of polynomials with simple roots, one of these roots, say $\lambda_{i,1}\in\bar{k}$,
$1\leq i\leq n$, being distinguished, and let \[
R_{i}\left(T\right)=\left(T-\lambda_{i,1}\right)^{-1}P_{i}\left(T\right)=\prod_{j=2}^{r_{i}}\left(T-\lambda_{i,j}\right),\quad i=1,\ldots,n.\]
 The nonsingular surface $S=S_{P_{1},\ldots,P_{n}}\subset\textrm{Spec}\left(\bar{k}\left[X_{0},X_{1},\ldots,X_{n+1}\right]\right)$
with equations \[
\left\{ \begin{array}{lcll}
X_{0}X_{j+1} & = & {\displaystyle {\displaystyle \prod_{i=1}^{j-1}}R_{i}\left(X_{i}\right)P_{j}\left(X_{j}\right)} & \textrm{for }1\leq j\leq n\\
\left(X_{j-1}-\lambda_{j-1,1}\right)X_{l+1} & = & X_{j}{\displaystyle \prod_{i=j}^{l-1}}R_{i}\left(X_{i}\right)P_{l}\left(X_{l}\right) & \textrm{for }2\leq j\leq l\leq n\end{array}\right.\]
 where $\prod_{i=j}^{l-1}R_{i}\left(X_{i}\right)=1$ if $j>l-1$,
is a Danielewski surface with base $\left(\bar{k}\left[x\right],x\right)$.
Indeed, $S\mid_{X_{*}}\simeq\textrm{Spec}\left(\bar{k}\left[X_{0},X_{0}^{-1}\right]\left[X_{1}\right]\right)$
whereas the irreducible components of the fiber $pr_{x}^{-1}\left(0\right)$
are the curves $C\left(\lambda_{1,1},\ldots,\lambda_{k,1},\ldots,\lambda_{m-1,j}\right)\simeq\ba_{\bar{k}}^{1}$
with equations \[
\left\{ X_{0}=0,\left(X_{i}=\lambda_{i,1}\right)_{i=1,\ldots,m-2},X_{m-1}=\lambda_{m-1,j}\right\} \quad m=1,\ldots,n,\, j=1,\ldots,r_{m}\]
 Moreover, the projection $p':S\rightarrow\textrm{Spec}\left(\bar{k}\left[X_{n+1}\right]\right)$
is a second $\ba^{1}$-fibration, and so $ML\left(S\right)=\bar{k}$.
In \cite{DFSemb}, we prove that a Danielewski surface with base $\left(\bar{k}\left[x\right],x\right)$
has a trivial Makar-Limanov invariant if and only if it is isomorphic
to a surface $S_{P_{1},\ldots,P_{n}}$. 
\end{example}
\begin{cor}
Given a morphism a Danielewski surfaces $\beta:S'\rightarrow S$ ,
we let \[
\beta:S'=S_{m}\stackrel{\beta_{m}}{\rightarrow}S_{m-1}\stackrel{\beta_{m-1}}{\rightarrow}\cdots\stackrel{\beta_{1}}{\rightarrow}S_{0}\stackrel{i}{\hookrightarrow}S\]
 be the factorization of $\beta$ into a sequence of fibered modifications
followed by an open immersion. Then the following hold.

a) If $\textrm{ML}\left(S_{k_{0}}\right)=k$ then this also holds
for $S_{k}$, $0\leq k\leq k_{0}-1$. 

b) If $\textrm{ML}\left(S\right)=k$ then this also holds for $S_{0}$. 
\end{cor}
\begin{proof}
The Makar-Limanov invariant of $S_{k_{0}}$ being trivial, we can
suppose that $S_{k_{0}}$ is isomorphic to the surface defined by
a $\left(\bar{k}\left[x\right],x\right)$-labelled comb $\gamma=\left(\Gamma,\sigma\right)$.
Then $S_{k_{0}-1}$ is isomorphic to the surface defined by a comb
obtained from $\gamma$ by blowing-down the leaves at the unique maximal
element $e$ of $\Gamma\setminus\textrm{Leaves}\left(\Gamma\right)$.
Thus $S_{k_{0}-1}$ has a trivial Makar-Limanov invariant by virtue
of (\ref{MLGDSThm}), and so, (a) follows by induction. To prove (b),
we can again suppose that $S\simeq S\left(\gamma\right)$ for a certain
$\left(\bar{k}\left[x\right],x\right)$-labelled comb $\gamma=\left(\Gamma,\sigma\right)$.
Then $S_{0}$ is isomorphic to the surface defined by the comb obtained
from $\gamma$ by deleting a subset $N\subset\textrm{Leaves}\left(\Gamma\right)$.
\end{proof}
\noindent As a consequence of the above description, we recover the
characterization of the Danielewski surfaces $S_{P,1}=\left\{ xz-P\left(y\right)=0\right\} $
in $\ba_{\bar{k}}^{3}=\textrm{Spec}\left(\bar{k}\left[x,y,z\right]\right)$
given in \cite{BML01}. 

\begin{cor}
\label{cor:ODS_carac} For a Danielewski surface $\pi:S\rightarrow\ba_{\bar{k}}^{1}$,
the following are equivalent.

a) $\textrm{ML}\left(S\right)=\bar{k}$ and the canonical sheaf $\omega_{S}$
of $S$ is trivial 

b) $S\simeq S\left(\gamma\right)$ for a suitable $\left(\bar{k}\left[x\right],x\right)$-labelled
comb $\gamma$ of height $h\leq1$. 

c) $S\simeq S_{P,1}\subset\ba_{\bar{k}}^{3}$ for a certain nonconstant
polynomial $P$ with simple roots. 
\end{cor}
\begin{proof}
By (\ref{EquivalentTreesTheorem}), we can suppose that $S=S\left(\gamma\right)$
for a certain essential $\left(\bar{k}\left[x\right],x\right)$-labelled
tree $\gamma$. The canonical sheaf $\omega_{S\left(\gamma\right)}$
is trivial if and only if if and only if the invertible sheaf $\mathcal{L}_{n,-\mu\left(\gamma\right)}$
on $X\left(n\right)$ is trivial (see the proof of (\ref{cor:TrivialCan_equiv_freeGa})).
By construction, this is the case if and only if all the leaves of
$\gamma$ are at the same level. On the other hand, by (\ref{thm:ML_general_carac}),
$\textrm{ML}\left(S\left(\gamma\right)\right)$ is trivial if and
only $\gamma$ is a comb. Thus (a) and (b) are equivalent since an
essential comb with all its leaves at the same level is either the
trivial tree or a nontrivial comb of height $1$. Finally, (b) and
(c) are equivalent by virtue of (\ref{exa:ODS}).
\end{proof}

\subsection{Cyclic quotients of Danielewski surfaces}

A closed fiber $S_{z_{0}}$ of an $\ba^{1}$-fibration $q:S\rightarrow Z$
on a normal affine surface is called \emph{degenerate} if it is not
isomorphic to $\ba_{k\left(z_{0}\right)}^{1}$. In \cite{Fie94},
Fieseler describes the structure of invariant neighbourhoods of degenerate
fibers of a quotient $\ba^{1}$-fibration $q:S\rightarrow Z=S//\mathbb{G}_{a,\bc}$
on a normal affine surface $S/\bc$ with a nontrivial $\mathbb{G}_{a,\bc}$-action.
This description, which remains valid for affine surfaces defined
over an algebraically closed field of caracteristic zero, can be reinterpreted
as follows. 

\begin{prop}
Let $q:S\rightarrow Z$ be an $\ba^{1}$-fibration on a normal affine
surface $S/\bar{k}$. If $S_{z_{0}}=q^{-1}\left(z_{0}\right)$ is
a degenerate fiber of $q$, then there exist a finite morphism $\phi:X\rightarrow Z'=\textrm{Spec}\left(\mathcal{O}_{Z,z_{0}}\right)$
totally ramified at the unique point $x_{0}$ over $z_{0}$ and a
DFS $\pi:\tilde{S}\rightarrow X$ with base $\left(X,x_{0}\right)$
equipped with an action of a cyclic group $\mathbb{Z}_{m}$ such that
$S\times_{Z}Z'=\tilde{S}/\mathbb{Z}_{m}$. 
\end{prop}
\begin{proof}
If $S_{z_{0}}$ is reduced then every irreducible component of $S_{z_{0}}$
is a connected component of $S_{z_{0}}$ and is isomorphic to $\ba_{\bar{k}}^{1}$
(see Lemma 1.2 in \cite{Fie94}). Thus $p_{2}:S'=S\times_{Z}Z'\rightarrow Z'$
is already a DFS with base $\left(Z',z_{0}\right)$. Otherwise, if
$S_{z_{0}}$ is irreducible but not reduced, say with multiplicity
$m\geq2$, then, by Theorem 1.7 in \cite{Fie94}, there exists a Galois
covering $\phi:X\rightarrow Z'$ of order $m$, \'etale over $Z'\setminus\left\{ z_{0}\right\} $
and totally ramified at the unique point $x_{0}$ over $z_{0}$, such
that the normalization $\tilde{S}$ of $\left(S'\times_{Z'}X\right)_{\textrm{red}}$
is a DFS $\pi:\tilde{S}\rightarrow X$, equipped with an action of
the Galois group $\mathbb{Z}_{m}$ of the covering $\phi$ such that
$S'\simeq\tilde{S}/\mathbb{Z}_{m}$. If $S_{z_{0}}$ is neither irreducible
nor reduced, with irreducible components $C_{i}$ of multiplicity
$m_{i}\geq1$, $i=1,\ldots,n$, then there exist DFS's $\pi_{i}:\tilde{S}_{i}\rightarrow X_{i}$
as above such that $S'_{i}=\left(S'\setminus S_{z_{0}}\right)\cup C_{i}$
is isomorphic to $\tilde{S}_{i}/\mathbb{Z}_{m_{i}}$ for every $i=1,\ldots,n$.
Letting $m$ be a common multiple of the $m_{i}$'s, these DFS's glue
in a similar way as in the proof of Theorem 1.8 in \cite{Fie94} to
a DFS $\pi:\tilde{S}\rightarrow Y$, equipped with an action of $\mathbb{Z}_{m}$
such that $S'\simeq\tilde{S}/\mathbb{Z}_{m}$. This completes the
proof. 
\end{proof}
\noindent The above description is local in nature. However, a similar
argument shows that a surface $S/\bar{k}$ equipped with an $\ba^{1}$-fibration
$q:S\rightarrow\ba_{\bar{k}}^{1}$ with at most one degenerate fiber
is isomorphic to cyclic quotient of a Danielewski surface. By Proposition
2.15 in \cite{Dub02}, a normal affine surface $S/\bar{k}$ with a
trivial Makar-Limanov admits an $\ba^{1}$-fibration $q:S\rightarrow\ba_{\bar{k}}^{1}$
with this property. So we obtain the following result. 

\begin{thm}
\label{thm:MLCarac}A normal affine surface $S/\bar{k}$ with a trivial
Makar-Limanov invariant is isomorphic to a cyclic quotient of a Danielewski
surface. 
\end{thm}
\noindent If the corresponding quotient morphism $\psi:S'\rightarrow S'/\mathbb{Z}_{m}\simeq S$
is \'etale in codimension $1$, then the Danielewski surface $S'$
has a trivial Makar-Limanov invariant too (see \emph{e.g.} \cite{Vas69})
. In general this is not the case, as shown by the following example. 

\begin{example}
A similar argument as in (\ref{exa:ODS}) shows that the Danielewski
surface with base $\left(\bc\left[x\right],x\right)$ defined by the
following labelled tree 

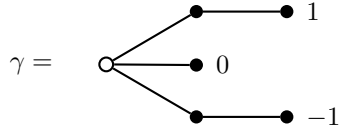
\begin{pspicture}(-2.2,-1)(8,1.2)

\rput(1.5,-0.02){$\gamma=$}

\rput(4,0){

\pstree[treemode=R,radius=2.5pt,treesep=0.5cm,levelsep=1.2cm]{\Tc{3pt}}{

\pstree{\TC*}{\TC*~[tnpos=r]{$1$}} \TC*~[tnpos=r]{$0$}\pstree{\TC*}{\TC*~[tnpos=r]{$-1$}}

}

}

\end{pspicture}

\noindent is isomorphic to the surface $pr_{x}:S'\rightarrow\textrm{Spec}\left(\bc\left[x\right]\right)$
in $\ba_{\bc}^{4}=\textrm{Spec}\left(\bc\left[x,y,z,u\right]\right)$
given by the equations \[
xz-y\left(y^{2}-1\right)=0,\quad\left(y^{2}-1\right)u-z^{2}=0,\quad xu-yz=0.\]
It comes equipped a $\mathbb{Z}_{2}$-action $\left(x,y,z,u\right)\mapsto\left(-x,-y,z,u\right)$.
The quotient $S/\mathbb{Z}_{m}$ is isomorphic to the nonsingular
surface $S'\subset\textrm{Spec}\left(\bc\left[a,b,c,d,e\right]\right)$
with equations \[
\left\{ \begin{array}{c}
ab=c^{2},\quad ad=\left(b-1\right)c,\quad ae=cd=b\left(b-1\right),\\
ce=bd,\quad\left(b-1\right)e=d^{2}\end{array}\right.,\]
 via the quotient morphism $\phi:S'\rightarrow S$ induced by the
ring homomorphism \begin{eqnarray*}
\bc\left[a,b,c,d,e\right] & \rightarrow & \bc\left[x,y,z,u\right],\left(a,b,c,d,e\right)\mapsto\left(x^{2},y^{2},xy,z,u\right).\end{eqnarray*}
 Note that $\phi$ is ramified along $pr_{x}^{-1}\left(0\right)$.
Since $\gamma$ is not comb, $\textrm{ML}\left(S'\right)$ is not
trivial. On the other hand, $S$ has trivial Makar-Limanov invariant
since it admits two $\ba^{1}$-fibrations $pr_{a}:S\rightarrow\textrm{Spec}\left(\bc\left[a\right]\right)$
and $pr_{e}:S\rightarrow\textrm{Spec}\left(\bc\left[e\right]\right)$
with distinct general fibers. 
\end{example}
\bibliographystyle{amsplain}

\end{document}